\documentclass[11pt]{amsart}
\usepackage{amsxtra}
\theoremstyle{plain}

\newcommand{\cleqn}{\setcounter{equation}{0}}
\newcommand{\clth}{\setcounter{theorem}{0}}
\newcommand {\sectionnew}[1]{\section{#1}\cleqn\clth}

\newtheorem{theorem}{Theorem}[section]
\newtheorem{lemma}[theorem]{Lemma}
\newtheorem{definition-theorem}[theorem]{Definition-Theorem}
\newtheorem{proposition}[theorem]{Proposition}
\newtheorem{corollary}[theorem]{Corollary}
\newtheorem{definition}[theorem]{Definition}
\newtheorem{example}[theorem]{Example}
\newtheorem{remark}[theorem]{Remark}
\newtheorem{notation}[theorem]{Notation}
\newtheorem{assumption}[theorem]{Assumption}
\newtheorem{lemma-definition}[theorem]{Lemma-Definition}
\newtheorem{lemma-notation}[theorem]{Lemma-Notation}
\newtheorem{question}[theorem]{Question}
\newcommand \bth[1] { \begin{theorem}\label{t#1} }
\newcommand \ble[1] { \begin{lemma}\label{l#1} }

\newcommand \bpr[1] { \begin{proposition}\label{p#1} }
\newcommand \bco[1] { \begin{corollary}\label{c#1} }
\newcommand \bde[1] { \begin{definition}\label{d#1}\rm }
\newcommand \bex[1] { \begin{example}\label{e#1}\rm }
\newcommand \bre[1] { \begin{remark}\label{r#1}\rm }

\newcommand \bnota[1] {\begin{notation}\label{n#1}\rm }
\newcommand \bas[1] { \begin{assumption}\label{a#1}\rm }

\newcommand \bln[1] { \begin{lemma-notation}\label{ln#1} }
\newcommand \bqu[1] { \begin{question}\label{q#1}\rm }
\newcommand {\eth} { \end{theorem} }
\newcommand {\ele} { \end{lemma} }

\newcommand {\epr} { \end{proposition} }
\newcommand {\eco} { \end{corollary} }
\newcommand {\ede} { \end{definition} }
\newcommand {\eex} { \end{example} }
\newcommand {\ere} { \end{remark} }
\newcommand {\enota} { \end{notation} }
\newcommand {\eas} {\end{assumption}}

\newcommand {\eln}{ \end{lemma-notation} }

\newcommand {\equ} {\end{question}}
\newcommand \thref[1]{Theorem \ref{t#1}}
\newcommand \leref[1]{Lemma \ref{l#1}}
\newcommand \prref[1]{Proposition \ref{p#1}}
\newcommand \coref[1]{Corollary \ref{c#1}}
\newcommand \deref[1]{Definition \ref{d#1}}
\newcommand \exref[1]{Example \ref{e#1}}

\newcommand \lb[1]{\label{#1}}
\newcommand \notaref[1]{Notation \ref{n#1}}

\newcommand \lnref[1]{Lemma-Notation \ref{ln#1}}

\def \Cset {{\mathbb C}}

\def \B  {{\mathcal{B}}}
\def \C  {{\mathcal{C}}}

\def \O  {{\mathcal{O}}}

\def \calP {{\mathcal{P}}}
\def \Pv {{\calP(1,v)}}
\def \I {{\mathcal I}}
\def \T {{\mathcal T}}

\def \Ga {\Gamma}

\def \fg{{\mathfrak g}}

\def \fk{{\mathfrak k}}

\def \fb{{\mathfrak b}}

\def \C{{\mathcal C}}

\def \hs {\hspace{.2in}}

\def \hs {\hspace{.2in}}

\def \sdot {\!\cdot\!}
\def \J {\mathcal J}

\def \th {\theta}
\def \bvs {({\bf v}, {\bf s})}
\def \bvsp {({\bf v}^\prime, {\bf s}^\prime)}

\def \calP {{\mathcal{P}}}
\def \Pv {{\calP(v_0,v)}}

\def \Wvv {W_{v_0}(v)}
\def \Yvv {\min(W_{v_0}(v))}
\def \Yvv {Y_{v_0}(v)}
\def \Zvv {Z_{v_0}(v)}

\def \bbvs { ({\bf v}_{v_0}({\bf s}), {\bf s})}

\def \Wvvp {W_{v_0}^{\prime}(v)}
\def \T {{\mathcal T}}

\begin{document}
\setlength{\baselineskip}{1.2\baselineskip}
\title[Some invariants of $K$-orbits in the flag variety]
{On some invariants of orbits in the flag variety under a symmetric subgroup}
\author[Sam Evens]{Sam Evens}
\address{Sam Evens, Department of Mathematics, University of Notre Dame, IL, USA}
\email{evens.1@nd.edu} 
\author[Jiang-Hua Lu]{Jiang-Hua Lu}
\address{Jiang-Hua Lu,
Department of Mathematics,
Hong Kong University,
Pokfulam Rd., Hong Kong}
\email{jhlu@maths.hku.hk}
\date{}
\begin{abstract} Let $G$ be a connected reductive algebraic group over an algebraically closed field ${\bf k}$ of 
characteristic not equal to $2$, let $\B$ be the
variety of all Borel subgroups of $G$, and let $K$ be a symmetric subgroup of $G$. Fixing a closed $K$-orbit in $\B$, 
we associate to every $K$-orbit on $\B$
some subsets of the Weyl group of $G$, and we study them as invariants of the $K$-orbits. 
When ${\bf k} = {\mathbb C}$, these invariants are used to determine 
when an orbit of a real form of $G$  and an orbit of a Borel subgroup of $G$ have non-empty intersection in $\B$.
We also characterize the invariants in terms of admissible paths in the set of $K$-orbits in $\B$.
\end{abstract}
\maketitle
\sectionnew{Introduction}\lb{intro}

Let $G$ be a connected reductive algebraic group over an algebraically closed field ${\bf k}$ of characteristic not equal to $2$, and let
$\B$ be the variety of Borel subgroups of $G$ with the conjugation action by $G$. Let
$\theta$ be an order $2$ automorphism of $G$.  The $\th$-fixed point subgroup $K = G^\theta$
acts on $\B$ with finitely
many orbits \cite{RS90, RS93, Springer85}. 

Let $V$ be the finite set of $K$-orbits in $\B$, and for $v \in V$, let $K(v) \subset \B$ be the corresponding $K$-orbit in $\B$.
For a subset $X$ of $\B$, let $\overline{X}$ be the Zariski closure of $X$ in $\B$.
The Bruhat order on $V$, denoted by $\leq$, is defined by 
\begin{equation}\lb{eq-Bruhat-V}
v_1 \leq v_2 \hs \mbox{if} \hs K(v_1) \subset \overline{K(v_2)}, \hs v_1, v_2 \in V.
\end{equation}
When ${\bf k} = \Cset$,  the geometry of the $K$-orbits and their closures in $\B$ 
plays an important role in the
representation theory of real forms of $G$ via the Beilinson-Bernstein
correspondence (see \cite{Mi93}). The poset $(V, \leq)$ and its application to representation theory have been the focus of
extensive studies (see \cite{Adams,  Adams-Barbasch-Vogan, Adams-Vogan, RS90, RS93, Springer85, Springer-invariants, Vogan-duke}).

Let $W$ be the canonical Weyl group of $G$ with the set $S$ of canonical generators, 
and let $M(W, S)$ be the corresponding monoid (see $\S$\ref{subsec-W-MW} and $\S$\ref{subsec-MWS}).
Among the structures on $V$ are the {\it monoidal action}  of $M(W, S)$ on $V$ (see $\S$\ref{subsec-MWS}), 
the {\it cross action} of $W$ on $V$ (see $\S$\ref{subsec-WV-cross}), 
and the map $\phi: V \to W$ defined by Springer (see $\S$\ref{subsec-Springer-map}). 
The programs from the {\it Atlas of Lie groups} (www.liegroups.org) facilitate
explicit computations in examples of
the cross action, the monoidal action, and the Springer map.

Let $V_0 = \{v \in V:  K(v) \subset \B \; \mbox{is closed}\}$. 
In this paper, we associate to every $v_0 \in V_0$ and $v \in V$ three subsets
$\Zvv \subset \Yvv \subset \Wvv$ of the Weyl group $W$. The set $\Wvv$ is defined both geometrically in terms of
intersections of $K(v)$ and Schubert varieties in $\B$ and combinatorially in terms of the Bruhat order on $V$ and the
monoidal action of $M(W, S)$ on $V$, and $\Yvv$ (resp. $\Zvv$) is defined to be the set of 
minimal (resp. minimal length) elements in $\Wvv$. 

One of our main results (see \thref{th-1-new} and \coref{co-Zvv-invariant}) says that 
for any $v_0 \in V_0$, the subset $\Yvv$ (or $\Zvv$) of $W$, together with
the Springer invariant $\phi(v) \in W$, completely determine $v \in V$. This result 
generalizes \cite[Theorem 5.2.2]{RS93} of Richardson-Springer in the case when $(G, \th)$ is of
Hermitian symmetric type (see $\S$\ref{subsec-hermitian}). 
When $V_0$ contains a single element $v_0$, the set $\Yvv$ for any $v \in V$ is a special case of an invariant
for $K(v)$ studied by Springer 
in \cite{Springer-invariants}.

When
$k=\Cset$ and $K$ is the complexification of a maximal compact subgroup of a real form $G_0$ of $G$, 
we use the Matsuki duality between $K$-orbits and $G_0$-orbits in $\B$
to determine when a 
$G_0$-orbit and a Bruhat cell in $\B$ have non-empty intersection.
The problem of determining when two such orbits intersect comes
 from Poisson geometry and  
served as the main motivation for this paper. 
See $\S$\ref{subsec-G0-B} and $\S$\ref{sec-G0-B}.

We introduce the notion of {\it admissible paths} from $v_0\in V_0$ to $v \in V$ by allowing only certain types of
cross and monoidal actions by simple generators of $W$. 
We then
characterize elements in $\Yvv$ (resp. $\Zvv$) as products of the simple generators in 
{\it minimal} (resp. {\it shortest}) admissible paths from $v_0$ to $v$.
See $\S$\ref{subsec-subexp-path} and  $\S$\ref{sec-Yvv-Zvv}.

The technical part of the paper is an analysis in
$\S$\ref{sec-analysis-Yvv} on the sets $\Yvv$ in relation to simple roots and their
types
relative to $v$.
Various examples are computed in $\S$\ref{sec-examples}.

Precise statements of the main results in the paper are given in $\S$\ref{sec-results}.

Most of the structures on $V$ are defined in \cite{RS90, RS93, Springer85} using a  
{\it standard pair}, i.e., a pair $(B, H)$, 
where $B \in \B$ and $H \subset B$ a maximal
torus of $G$ such that $\th(B) = B$ and $\th(H) = H$.
However, to formulate the sets $\Zvv \subset \Yvv\subset \Wvv$,
 it is crucial that these structures can be defined  
independently of the standard pairs. 
Following \cite[$\S$1.7]{RS93},  
we review in $\S$\ref{sec-review-I} 
and $\S$\ref{sec-review-II}
the canonical definitions of 
all the structures on $V$ needed in the paper. 

\subsection*{Acknowledgments} 
We would like to thank Xuhua He and George McNinch for helpful discussions.
We thank the organizers
of {\it the Atlas of Lie Groups} meeting in Utah in July 2010, which
helped the  development of this paper, and also thank Jeff Adams, Peter
Trapa, and David Vogan for very useful discussions.  
We are especially grateful
to Scott Crofts, who wrote a program allowing us to compute the invariants 
$\Zvv$ and $\Yvv$ explicitly in examples. The work described in this paper by the
first author was partially supported by NSA grants H98230-08-0023 and
 H98230-11-1-0151, and that
by the second author  by the RGC grants HKU 7034/05P
and 7037/07P of the Hong Kong SAR, China. 

\sectionnew{Statements of Results}\lb{sec-results}

\subsection{Geometric and combinatorial definitions of $\Wvv\subset W$}\lb{subsec-twodef-Wvv}

For $v_0 \in V_0$ and $v\in V$,  
we define $\Wvv \subset W$ by
\begin{equation}\lb{eq-Wvv-intro}
\Wvv = \{w \in W: K(v) \cap \overline{B(w)} \neq \emptyset\},
\end{equation}
where $B$ is any Borel subgroup of $G$ contained in $K(v_0)$, and for $w \in W$, $B(w)$ is the corresponding $B$-orbit in $\B$
(see 
$\S$\ref{subsec-B-orbits}). The set $\Wvv$ depends only on $v_0$ and $v$ and not on the choice of $B \in K(v_0)$ (see \leref{le-Wvv-well}).

Denote the monoidal action of $M(W, S)=\{m(w): w \in W\}$ on $V$ by  $m(w)\sdot v$ for $w \in W$ and $v \in V$
(see $\S$\ref{subsec-MWS}). Our \leref{le-Wvv-combinatorial} gives the following combinatorial description of the
set $\Wvv$ for $v_0 \in V$ and $v \in V$:
\begin{equation}\lb{eq-Wvv-intro-2}
\Wvv = \{w \in W: \; v \leq m(w)\sdot v_0\} \subset W,
\end{equation}
where recall that $\leq$ is the Bruhat order on $V$ defined in \eqref{eq-Bruhat-V}.

\subsection{The subsets $\Zvv \subset \Yvv$ of $\Wvv$}\lb{subsec-ZY}
Let $l: W \to {\mathbb N}$ and $\leq$ be respectively the length function and the Bruhat order  
on $W$ as a Coxeter group (see $\S$\ref{subsec-W-MW} and $\S$\ref{subsec-Bruhat-W}). 
Let  $W_1$ be a subset of $W$. An element $w \in W_1$ is said to be {\it minimal}
if for any $w_1 \in W_1$, $w_1 \leq w$ implies that $w_1 = w$.  
The set of all minimal elements in $W_1$ will be denoted by $\min(W_1)$. An element $w \in W_1$ is said to have {\it minimal length}
if $l(w) \leq l(w_1)$ for every $w_1 \in W_1$. The set of all minimal length elements in $W_1$ is denoted by $\min_l(W_1)$. 

For $v_0 \in V_0$ and $v \in V$, define
\begin{align}
\lb{eq-Yvv-intro}
\Yvv &= \min (\Wvv) \subset \Wvv,\\
\lb{eq-Zvv-intro}
\Zvv &= {\rm min}_l(\Wvv) ={\rm min}_l(\Yvv) \subset \Yvv.
\end{align} 
We prove in \leref{le-Wvv-Yvv} that for any $v_0 \in V$ and $v \in V$,   
$\Wvv$ is determined by its subset $\Yvv$ in the sense that for any $w \in W$, 
\begin{equation}\lb{eq-w-y}
w \in \Wvv \hs \mbox{iff} \hs  y \leq w \;\; \mbox{for some} \;y \in \Yvv.
\end{equation}  
Thus both $\Zvv$ and $\Wvv$ are determined by $\Yvv$. The example in $\S$\ref{subsec-an-ex} shows that $\Zvv$ may be 
a proper subset of $\Yvv$. 

Using the programs available at the {\it Atlas of Lie groups}, Scott Crofts has
written a  program that allows one to compute the sets $\Yvv$ 
explicitly in examples.

\subsection{Intersections of real group orbits and Bruhat cells in $\B$}\lb{subsec-G0-B}
Assume now that ${\bf k} = {\mathbb C}$ and that $K$ is the complexification of a maximal compact
subgroup of a real form $G_0$ of $G$. For $v \in V$, let $G_0(v) \subset \B$ be the $G_0$-orbit in $\B$ that is
dual to 
$K(v) \subset \B$ under the Matsuki duality \cite{Uzawa-flow}  between $G_0$-orbits and $K$-orbits in $\B$.
The following \thref{th-2-new}  is the
first main result of this paper. 

\bth{th-2-new}
Let $v_0 \in V_0$ and let $B \in K(v_0)$. Then for  $v \in V$ and $w \in W$,
\[
G_0(v) \cap B(w)  \neq \emptyset \hs \mbox{iff} \hs w \in \Wvv.
\]
\eth

Our motivation for \thref{th-2-new} 
comes from Poisson geometry: when $G_0$ is connected, it is
shown in \cite{Foth-Lu:real} that there is a Poisson structure $\Pi$ on $\B$ such that
the connected components of intersections of $G_0$-orbits and $B$-orbits in $\B$ are
precisely the $H_0$-orbits of symplectic leaves of $\Pi$, where $H_0=B \cap G_0$ is a maximally
compact Cartan subgroup of $G_0$. Thus, one first needs to know when the intersection of a $G_0$-orbit and a $B$-orbit in $\B$ 
is non-empty. In view of \eqref{eq-w-y},
\thref{th-2-new} gives a complete answer to this question in terms of the set $\Yvv$ and the Bruhat order on $W$.
Further applications to the Poisson structure $\Pi$ on $\B$ will appear elsewhere.

\subsection{$\Yvv$ and $\Zvv$ as invariants of $K(v)$}\lb{subsec-as-invariants}
Our second main result, the following \thref{th-1-new}, implies that
the pair $(\phi(v), \Yvv)$ forms a complete invariant for $v \in V$, where $\phi: V \to W$ is the Springer map
(see $\S$\ref{subsec-Springer-map}).

\bth{th-1-new} Let $v_0 \in V_0$ and $v, v^\prime \in V$. If $\phi(v) = \phi(v^\prime)$ and
$\Yvv\cap Y_{v_0}(v^\prime) \neq\emptyset,$ then $v = v^\prime$.
\eth

Since $\Zvv \subset \Yvv$ for any $v_0 \in V_0$ and $v \in V$, one also has

\bco{co-Zvv-invariant}
Let $v_0 \in V_0$ and $v, v^\prime \in V$. If $\phi(v) = \phi(v^\prime)$ and
$\Zvv\cap Z_{v_0}(v^\prime) \neq\emptyset,$ then $v = v^\prime$.
\eco

\subsection{Subexpressions and admissible paths}\lb{subsec-subexp-path}
Reduced decompositions for elements in $V$ and their subexpressions are introduced in \cite[$\S$5]{RS90} and
\cite[$\S$3, $\S$4]{RS93} (and see $\S$\ref{subsec-reduced} for details). In particular, if $v, v^\prime \in V$ and if  
$\bvsp$ is any reduced decomposition of $v^\prime$, then \cite[Proposition 4.4]{RS93}  $v \leq v^\prime$ if and only if
there exists a subexpression of $\bvsp$ with final term $v$. 

Let $v_0 \in V_0$ and $v \in V$. We show  that if $y \in \Yvv$, then every reduced decomposition
of $m(y)\sdot v_0$ coming from a reduced word of $y$ has 
{\it exactly one} subexpression with final term $v$. See \prref{pr-unique} for details. 

For $v_0 \in V_0$ and $v \in V$, we define an {\it admissible path} from $v_0$ to $v$ to be a pair 
$({\bf v}, {\bf s})$, where ${\bf v} = (v_0, v_1, \ldots, v_k)$ is a sequence in $V$ and ${\bf s} = (s_{1}, \ldots, s_{k})$
is a sequence in $S$, such that for each $j \in [1,k]$, $v_j$ is a certain type of either the cross action or the monoidal action of
$s_j$ on $v_{j-1}$ (see \deref{de-1-path}). 
For such a path $({\bf v}, {\bf s})$, let $y({\bf v}, {\bf s}) =s_k \cdots s_1 \in W$.
We also define the set ${\mathcal P}_{\rm min}(v_0, v)$ (resp. ${\mathcal P}_{{\rm short}}(v_0, v)$)
 of  {\it minimal} (resp. {\it shortest}) paths from  $v_0$ to $v$. 
We prove (see \coref{co-Yvv} and \coref{co-Zvv}) that for any $v_0 \in V_0$ and $v \in V$, 
\begin{align*}
\Yvv &= \{y\bvs: \bvs \in {\mathcal P}_{\rm min}(v_0, v)\}\\
\Zvv &= \{y\bvs: \bvs \in {\mathcal P}_{\rm short}(v_0, v)\}.
\end{align*}
We believe that the minimal and shortest paths defined in this paper will have other applications to the study of $K$-orbit closures in $\B$.

\sectionnew{Review on $K$-orbits in $\B$, I}\lb{sec-review-I}
Following \cite[$\S$1.7]{RS93}, we review in this section the canonical Weyl group $W$ of $G$ with its set $S$ of canonical generators
and the  action of the monoid  $M(W, S)$
on $V$. More structures on $V$ will be reviewed in $\S$\ref{sec-review-II}.

\subsection{Notation} 
If $Q$ is a subgroup of $G$ and $g \in G$, $Q^g$ denotes the subgroup $g^{-1}Qg$ of $G$. We will consider the right action
of $G$ on the flag variety $\B$ by $\B \times G \to \B: (B, g) \mapsto B^g$ for $B \in \B$ and $g \in G$.

If a group $L$ acts on a set $X$ from the left (resp. right), we will denote by $L \backslash X$ (resp. $X/L$) the set 
of $L$-orbits on $X$. 

Throughout the paper and unless otherwise specified,
for a subset $X$ of $G$ (resp. $\B$ and $\B \times \B$), $\overline{X}$ denotes the Zariski closure of $X$ in $G$ (resp. 
$\B$ and
$\B \times \B$).
\subsection{The variety $\C$}\lb{subsec-C} 
Let $\C$ be the set of all pairs $(B, H)$, where $B \in \B$ and $H \subset B$ is a maximal torus of $G$.
Then $G$ acts on $\C$ transitively from the right by
\begin{equation}\lb{eq-G-C}
\C \times G \longrightarrow \C: \; \; (B, H)^g = (B^g, \; H^g), \hs (B, H) \in \C, \, g \in G.
\end{equation}
The stabilizer subgroup of $G$ at $(B, H) \in \C$ is
$H$. Thus for each $(B, H) \in \C$, one has the $G$-equivariant identification
\begin{equation}\lb{eq-CBH}
C_{B, H}: \; \; \; H\backslash G \longrightarrow \C: \; \; \; Hg \longmapsto (B^g, H^g), \hs g \in G.
\end{equation}

For $(B, H) \in \C$, let $N_G(H)$ be the normalizer of $H$ in $G$, let
$W_H = N_G(H)/H$ be the Weyl group of $(G, H)$, and let 
$S_{B, H}$ be the set of generators of $W_H$ defined by the simple roots of $B$ relative to $H$.  

Let $(B, H), (B^\prime, H^\prime) \in \C$. Let $g \in G$  be such that 
$B^\prime = B^{g}$ and $H^\prime = H^{g}$. Since $g$ is unique up to the left multiplication by 
an element in $H$, one  has a well-defined isomorphism of tori
\begin{equation}\lb{eq-T-BBHH}
T_{B, H}^{B^\prime, H^\prime}: \; \; H \longrightarrow H^\prime: \; \; h \longmapsto g^{-1} h g, \hs h \in H.
\end{equation}
Moreover, although the map $N_G(H) \to N_G(H^\prime): n \to g^{-1} n g$ depends on the choice of $g$, 
the  group isomorphism
\[
\eta_{B, H}^{B^\prime, H^\prime}: \; \; W_H \longrightarrow W_{H^\prime}: \; \; nH \longmapsto g^{-1} n g H^\prime, \hs n \in N_G(H)
\]
does not and is thus well-defined. Since $\eta_{B, H}^{B^\prime, H^\prime}(S_{B, H}) = S_{B^\prime, H^\prime}$, one has the isomorphism
$\eta_{B, H}^{B^\prime, H^\prime}: (W_H, S_{B, H}) \to (W_{H^\prime}, S_{B^\prime, H^\prime})$ of
Coxeter groups.

\subsection{The canonical Weyl group}\lb{subsec-W-MW} 
Let $W=(\B \times \B)/G$ be the set of $G$-orbits on $\B \times \B$ for the diagonal $G$-action, and let 
$p: \B \times \B \to W$ be the natural projection. 
Let $(B, H) \in \C$. Then the map
\begin{equation}\lb{eq-etaBH}
\eta_{B, H}: \; \; \; W_H \longrightarrow W: \; \; nH \longmapsto p(B^n, B), \hs n \in N_G(H)
\end{equation}
is bijective.
It is straightforward to check that for another $(B^\prime, H^\prime) \in \C$, 
\[
\eta_{B^\prime, H^\prime}^{-1} \circ \eta_{B, H}=\eta_{B, H}^{B^\prime, H^\prime}: \;\;\; W_H \longrightarrow W_{H^\prime}.
\]
Thus there
is a well-defined group structure on $W$ such that $\eta_{B, H}: W_H \to W$ is 
a group isomorphism for every $(B, H) \in \C$. Let $S = \eta_{B, H}(S_{B, H}) \subset W$ for any $(B, H) \in \C$. Then
$S$ is a set of generators for $W$ and is independent of the choice of $(B, H) \in \C$.
The Coxeter group $(W, S)$  is called 
the {\it canonical Weyl group} of $G$.

For $w \in W$, 
a reduced word of $w$ is a shortest expression $w = s_1 \cdots s_{l(w)}$ of $w$ as a product of
elements in $S$, and $l(w)$ is called the length of $w$. 

\subsection{The Bruhat order on $W$}\lb{subsec-Bruhat-W}
For $w \in W$, let 
\[
\O(w) =\{(B^\prime, B) \in \B \times \B: \, p(B^\prime, B) = w\} \subset \B \times \B
\]
 be the corresponding $G$-orbit in $\B \times \B$. 
The Bruhat order on $W$ 
is defined by
\[
w \leq w^\prime \hs \mbox{if} \hs \O(w) \subset \overline{\O(w^\prime)}, \; \; \; w, w^\prime \in W.
\]
If $w \leq w^\prime$ and $w \neq w^\prime$, we will write $w < w^\prime$ or $w^\prime > w$.
It is well-known \cite[Th\'eor\`eme 3.13]{Borel-Tits} that for $w, w^\prime \in W$,  $w^\prime \leq w$ if and only if
there is a reduced word $w = s_1 \cdots s_{l(w)}$ of $w$ such that $w^\prime = s_{i_1} \cdots s_{i_p}$ for some
$1 \leq i_1 < \cdots < i_p \leq l(w)$. 

For a subset $W_1$ of $W$, recall from $\S$\ref{subsec-ZY} that $\min(W_1)$ is the set of minimal 
elements
in $W_1$ with respect to the Bruhat order and that $\min_l(W_1)$ is the set
of minimal length elements in $W_1$.
The following \leref{le-useful} will be used in the proof of \prref{pr-Yvv-Wvvp}.
 
\ble{le-useful}
Let $W_1$ and $W_2$ be two subsets of $W$ such that $\min (W_2) \subset W_1 \subset W_2$. 
Then $\min(W_1)= \min (W_2)$ and $\min_l(W_1)= \min_l(W_2)$.
\ele

\subsection{Orbits in $\B$ under a Borel subgroup}\lb{subsec-B-orbits}
Let $B \in \B$. The set of $B$-orbits in $\B$ is naturally indexed by
$W$. Indeed, for $w \in W$, let 
\begin{equation}\lb{eq-Bw}
B(w) =\{B^\prime \in \B: \, p(B^\prime, B) = w\} \subset \B.
\end{equation}
Then $B(w)$ is a single $B$-orbit in $\B$, and the map $w \mapsto B(w)$  is a bijection from $W$ to the set of all $B$-orbits in $\B$.

Define $q_B: G \to \B$ by $q_B(g) = B^g$ for $g \in G$. For $w \in W$, let 
\begin{equation}\lb{eq-BwB}
BwB = q_B^{-1}(B(w)) =\{g \in G: p(B^g, B) = w\} \subset G.
\end{equation}
Then $BwB$ be the single $(B, B)$-double coset in $G$. Moreover, 
for $w, w^\prime \in W$,
\[
w \leq w^\prime \hs \mbox{iff} \hs B(w) \subset \overline{B(w^\prime)} \hs \mbox{iff} \hs BwB \subset \overline{Bw^\prime B}.
\]

\subsection{The monoidal action of $M(W, S)$ on $V$}\lb{subsec-MWS} Our references for this subsection are
\cite[$\S$4]{RS90} and \cite[$\S$3]{RS93}.
The monoid $M(W, S)$ associated to the Coxeter group $(W, S)$ is  
$M(W, S) = \{m(w): w \in W\}$
as a set,
with the monoidal product given by
\begin{equation}\lb{eq-msmw}
m(s) m(w) = \begin{cases} m(s w), & \;\;\; \mbox{if} \;\; s w > w\\ m(w), & \;\;\; \mbox{if} \;\; s w < w\end{cases}, \hs s \in S,\; w \in W.
\end{equation}
Alternatively, let $B \in \B$, and let $BwB = q_B^{-1}(B(w))\subset G$ for $w \in W$ as in \eqref{eq-BwB}. Define $\ast: W \times W \to W$ such that 
\begin{equation}\lb{eq-BwwB}
\overline{B (w_1 \ast w_2) B} = \overline{(Bw_1B)(Bw_2B)}, \hs w_1, w_2 \in W.
\end{equation}
Then  $\ast$ is a monoidal product on $W$, independent of
the choice of $B \in \B$, and $(W, \ast) \to M(W, S): w \mapsto m(w)$ is an isomorphism of monoids
\cite[Proposition 3.18]{Borel-Tits}.

Let $q_K: G \to G/K$ be the natural projection.
Let $B \in \B$, and for $v \in V$, let 
\[
BvK = q_B^{-1}(K(v)) \subset G.
\]
Then $q_K(BvK)$, being a $B$-orbit in $G/K$, is irreducible for each $v \in V$.
For $w \in W$ and $v \in V$, since $q_K((BwB) (BvK)) \subset G/K$ is irreducible  and is a finite union
of $B$-orbits in $G/K$,  there is a unique element in $V$, denoted by $m(w)\sdot v$, such that
$\overline{q_K((BwB) (BvK)} = \overline{q_K(B(m(w)\sdot v)K)}$, where $\,\bar{}\,$ denotes the Zariski closure in $G/K$. Consequently, one has,
for any $w \in W$ and $v \in V$,
\begin{equation}\lb{eq-BK-monoidal}
\overline{B (m(w)\sdot v) K} = \overline{(BwB)(BvK)}.
\end{equation}

\ble{le-monoidal-action-well} \cite[$\S$4]{RS90}
The map
\begin{equation}\lb{eq-MWS-V}
M(W, S) \times V \longrightarrow V: \; \; \; (m(w), v) \longmapsto m(w) \sdot v, \hs w \in W, v \in V,
\end{equation}
is independent of the choice of $B \in \B$ and defines a monoidal action of $M(W, S)$ on $V$.
\ele

\begin{proof} Let $B \in \B, g_0 \in G$, and $B^\prime = B^{g_0}$. It follows
from definitions that for any $w \in W$ and $v \in V$,
$B^\prime w B^\prime = g_0^{-1} BwB g_0$ and $B^\prime v K = g_0^{-1} BvK$. 
Thus, the map in \eqref{eq-MWS-V} is independent of the choice of $B$. Choosing any $B \in \B$ and 
using the isomorphism $(W, \ast) \to M(W, S)$ of monoids, one sees from \eqref{eq-BwwB} and \eqref{eq-BK-monoidal} that
the map in \eqref{eq-MWS-V} defines an action of $M(W, S)$ on $V$.
\end{proof}

\ble{le-ast-bar}
For any $B \in \B, v \in V$, and $w \in W$, one has
\[
\overline{B (m(w)\sdot v)K}=\overline{BwB} \;\;\overline{BvK}.
\]
\ele

\begin{proof} Let $w=s_1 s_2\cdots s_l$ be a reduced word, and for each $1 \leq j \leq l$, let $P_j = 
B \cup Bs_jB$ so that $P_j$ is a parabolic subgroup of $G$. Then $\overline{BwB} = P_1P_2\cdots P_j$
(see, for example \cite[Th\'eor\`eme 3.13]{Borel-Tits}). By repeatedly applying
\cite[Lemma 2.12]{Borel-Tits}, one sees that $\overline{BwB} \;\;\overline{BvK}$ is closed in $G$.
Thus $\overline{BwB} \;\;\overline{BvK} = \overline{(BwB)(BvK)} = \overline{B (m(w)\sdot v)K}$.
\end{proof}

Recall that $V_0=\{v_0 \in V: K(v_0) \subset \B \; \mbox{is closed}\}$. 
For a subset $X$ of $\B$, let
\[
X\cdot K = \{B^k: \; B \in X, k \in K\}.
\]

\ble{le-BwK} 
Let $v_0 \in V_0$ and let $B \in K(v_0)$. Then for any $w \in W$,
\[
\overline{K(m(w)\sdot v_0)} = \overline{B(w) \cdot K}=\overline{B(w)}\cdot K.
\]
\ele

\begin{proof} Since $K(v_0)$ is closed in $\B$, the double coset 
$Bv_0 K = BK = q_B^{-1}(K(v_0))$
is closed in $G$. By \leref{le-ast-bar}, one has
\[
\overline{B(m(w)\sdot v_0)K} =\overline{B w B} \; \overline{BK} =\overline{BwB} \; K.
\]
Applying $q_B: G \to \B$, one proves \leref{le-BwK}.
\end{proof}

\ble{le-y1-y}
If $w, w^\prime \in W$ and $v, v^\prime \in V$ are such that $w \leq w^\prime$, $v \leq v^\prime$, then
\[
m(w) \sdot v \leq m(w^\prime) \sdot v \hs \mbox{and} \hs
m(w)\sdot v \leq m(w) \sdot v^\prime.
\]
\ele

\begin{proof} This is immediate from \leref{le-ast-bar}.
\end{proof}

\sectionnew{The geometrical and combinatorial definitions  of $\Wvv$}\lb{sec-geom}

\subsection{The two definitions  of $\Wvv$}\lb{subsec-geom-Wvv}
Recall from  $\S$\ref{subsec-twodef-Wvv} that for $v_0 \in V_0$ and $v \in V$, the subset $\Wvv$ of $W$ is defined by
\[
\Wvv = \{w \in W: K(v) \cap \overline{B(w)} \neq \emptyset\},
\]
where $B$ is any Borel subgroup of $G$ contained in $K(v_0)$.

\ble{le-Wvv-well}
For any $v_0 \in V$ and $v \in V$, the set $\Wvv$ is independent of the choice of $B \in K(v_0)$.
\ele

\begin{proof} Let $B, B^\prime \in K(v_0)$ and let $B^\prime = B^k$ for some $k \in K$. Let $w \in W$. Then 
$B^\prime (w)=\{B_1^k: B_1 \in B(w)\}$. Thus $K(v) \cap \overline{B^\prime(w)} \neq \emptyset$ if and only if $K(v) \cap 
\overline{B(w)} \neq \emptyset$.
\end{proof}

We now have the following combinatorial interpretation of $\Wvv$.

\ble{le-Wvv-combinatorial}
For any $v_0 \in V_0$ and $v \in V$, one has
\[
\Wvv = \{w \in W: \; v \leq m(w)\sdot v_0\}.
\]
\ele

\begin{proof} Let $w \in W$.
By \leref{le-BwK} and by the  definition of the Bruhat order on $V$,  $v \leq m(w) \sdot v_0$ if and only if
$K(v) \subset \overline{K(m(w)\sdot v_0)} = \overline{B(w)} \cdot K$,
which is  equivalent to $K(v) \cap \overline{B(w)} \neq \emptyset$.
\end{proof}

\subsection{Properties of $\Wvv$} 
Recall from $\S$\ref{subsec-ZY} that for $v_0 \in V_0$ and $v \in V$, $\Yvv$ is the set of minimal elements in $\Wvv$ with
respect to the Bruhat order on $W$.

\ble{le-Wvv-Yvv} Let $v_0 \in V_0$, $v \in V$, and $w \in W$.
Then $w \in \Wvv$ if and only if $w \geq y$ for some $y \in \Yvv$.
\ele

\begin{proof}
It follows from the definition of $\Yvv$ that $w \geq y$ for some $y \in \Yvv$. Conversely,
assume  that $w \geq y$ for some $y \in \Yvv$. Then $K(v) \cap \overline{B(y)} \neq \emptyset$ and 
$\overline{B(y)} \subset \overline{B(w)}$. Thus $K(v) \cap \overline{B(w)} \neq\emptyset$.
\end{proof} 



\ble{le-Wv-W}
For $v_0 \in V_0$ and $v \in V$, one has $\Wvv = W$ if and only if $v = v_0$.
\ele

\begin{proof} Let $1$ be the identity element of $W$. 
It is clear that $1 \in W_{v_0}(v_0)$, so $W_{v_0}(v_0)=W$ by \leref{le-Wvv-Yvv}.
  Assume that $v \in V$ is such that $\Wvv = W$. Then $1 \in \Wvv$, so $v \leq v_0$. Since
$K(v_0) \subset \B$ is closed, we have $v = v_0$.
\end{proof}

\subsection{The set $\Wvvp$}\lb{subsec-Wvvp} Fix $v_0 \in V_0$ and let $B \in K(v_0)$. For $v \in V$, let
\begin{equation}\lb{eq-Wvvp}
\Wvvp =\{w \in W: \; K(v) \cap B(w) \neq \emptyset\}.
\end{equation}
By the proof of \leref{le-Wvv-well},  $\Wvvp$ is independent of the choice of $B$ in  $K(v_0)$. One also sees from the definition of 
$B(w)$ in \eqref{eq-Bw} that for any $v_0 \in V_0$ and $v \in V$,
\begin{equation}\lb{eq-Wvvp-1}
\Wvvp = \{w \in W: \; w = p(B^\prime, B) \; \mbox{for some} \; B^\prime \in K(v), \, B \in K(v_0)\}.
\end{equation}
The following \prref{pr-Yvv-Wvvp} expresses $\Yvv$ and $\Zvv$ in terms of $\Wvvp$.

\bpr{pr-Yvv-Wvvp}
For any $v_0 \in V_0$ and $v \in V$, one has
\begin{equation}\lb{eq-Yvv-Wvvp}
\Yvv \subset \Wvvp \subset \Wvv.
\end{equation}
Consequently, $\Yvv = \min (\Wvvp)$ and $\Zvv = \min_l(\Wvvp)$.
\epr

\begin{proof} Clearly $\Wvvp \subset \Wvv$. 
Let $w \in \Yvv$. Then 
$K(v) \cap \overline{B(w)} \neq \emptyset$ and  $K(v) \cap B(w^\prime) =\emptyset$ for any $w^\prime\in W$ 
such that $w^\prime < w$. Thus $K(v) \cap B(w) \neq \emptyset$ and $w \in \Wvvp$. 
By \leref{le-useful}, $\Yvv = \min (\Wvvp)$, and $\Zvv = \min_l(\Wvvp)$.
\end{proof}

\bex{ex-G-diagonal} Let $\tilde{G} = G \times G$ and let 
\[
\tilde{\th}: \; \; \tilde{G} \longrightarrow \tilde{G}: \; \; \; \tilde{\th}(g_1, g_2) = (g_2, g_1), \hs (g_1, g_2) \in \tilde{G}.
\]
Then the fixed point subgroup $\tilde{K}$ of $\tilde{\th}$ in $\tilde{G}$ is $\tilde{K} = \{(g, g): g \in G\}$, so
the set $\tilde{V}$ of $\tilde{K}$-orbits in $\tilde{\B} = \B \times \B$ is $W$.  
For $w \in W$, let 
$\tilde{K}(w) = \O(w)$, where we recall that $\O(w)$ is the $G$-orbit in $\tilde{\B}$ for the diagonal action.  Then
the only $\tilde{v}_0 \in \tilde{V}= W$ such that $\tilde{K}(\tilde{v}_0)$ is closed in $\tilde{\B}$ is  $\tilde{v}_0=1$, the identity element
of $W$. 
Let $\tilde{W} = W \times W$, and for  $w \in W=\tilde{V}$, 
let $\tilde{W}_{\tilde{v}_0}^\prime(w)\subset \tilde{W}_{\tilde{v}_0}(w) \subset \tilde{W}$ be defined as in \eqref{eq-Wvv-intro-2}
and \eqref{eq-Wvvp} but for the pair $(\tilde{G}, \tilde{\th})$. Then it is easy to see that, for any $w \in W=\tilde{V}$, 
\[
\tilde{W}_{\tilde{v}_0}(w)=\{(w_1, w_2) \in W \times W: \; w \leq w_1 \ast w_2^{-1}\},
\]
where recall from $\S$\ref{subsec-MWS} that $\ast$ is the monoidal product on $W$, while
\[
\tilde{W}_{\tilde{v}_0}^\prime(w)=\{(w_1, w_2) \in W \times W: \; BwB \subset Bw_1Bw_2^{-1}B\},
\]
where $B$ is any Borel subgroup of $G$.
The set $\tilde{W}_{\tilde{v}_0}^\prime(w)$ can be computed by choosing any reduced word  for $w_1$
and using inductively the fact that, for $s \in S$ and $u \in W$, $BsBuB = BsuB$ if $su > u$ and $BsBuB = BsuB \cup BuB$ if $su < u$.
See \cite[Remarques 3.19]{Borel-Tits}. Moreover, for any $w \in W$, 
\begin{align*}
{\rm min} (\tilde{W}_{\tilde{v}_0}(w)) &= {\rm min}(\tilde{W}_{\tilde{v}_0}^\prime(w)) \\
&= \{(w_1, w_2) \in W \times W: \; 
w = w_1w_2^{-1}, \; l(w) = l(w_1) + l(w_2)\},
\end{align*}
and all elements in ${\rm min} (\tilde{W}_{\tilde{v}_0}(w))$ have the same length, namely $l(w)$.
\eex

\subsection{The element $m(w)\sdot v_0$} 
Fix $v_0 \in V_0$ and let $B \in K(v_0)$. 
For $w \in W$, let
\begin{equation}\lb{eq-WV0}
V_w=\{v \in V: \; K(v) \cap B(w) \neq \emptyset\}.
\end{equation}
Then $V_w$ depends only on $w$ and not
on the choice of $B \in K(v_0)$.

\bpr{pr-mwv0}
Let $v_0 \in V_0$, and let $B \in K(v_0)$. Then for every $w\in W$, 

1) $B(w) \cap K(m(w)\sdot v_0)$ is dense in 
$B(w)$,

2) $m(w)\sdot v_0$ is the unique maximal element in $V_w$ with respect to the Bruhat order on $V$.
\epr

\begin{proof} 
Since $B(w)$ is irreducible and since $\B = \bigsqcup_{v \in V} K(v)$,
there exists a unique $v_w \in V$ such that $B(w) \cap K(v_w)$ is dense in $B(w)$. Moreover, if
$v \in V_w$, then 
\[
\emptyset \neq B(w) \cap K(v) \subset B(w) \subset \overline{B(w)} = \overline{B(w) \cap K(v_w)}\subset \overline{K(v_w)},
\]
and so $v \leq v_w$. This shows that $v_w$ is the unique maximal element in $V_w$ with respect to the Bruhat 
order on $V$. 
It remains to show that $v_w = m(w)\sdot v_0$. 

Since $w \in  W_{v_0}(v_w)$, $v_w \leq m(w)\sdot v_0$. By \leref{le-BwK}, 
$\overline{K(m(w)\sdot v_0)} = \overline{B(w)\cdot K}$, 
so $K(m(w)\sdot v_0) \cap (B(w)\cdot K) \neq \emptyset$, and hence $K(m(w)\sdot v_0) \cap B(w) \neq \emptyset$.
Thus $m(w)\sdot v_0 \in V_w$ and $m(w)\sdot v_0 \leq v_w$. Hence $v_w=m(w)\sdot v_0$.
\end{proof}

\sectionnew{Intersections of $G_0$-orbits and $B$-orbits}\lb{sec-G0-B}
In this section, we assume that ${\bf k} = {\mathbb C}$ and that $G_0 =\{g \in G: \tau(g) = g\}$ 
is a real form of $G$, where $\tau$ is 
an anti-holomorphic involution on $G$. Let $\sigma$ be a Cartan involution of $G$  commuting with $\tau$, and let
$\theta = \sigma\tau$. Let $K = G^\theta$ and $K_0 = G_0 \cap K$. Then  $K_0$ is a maximal compact subgroup of $G_0$ and $K$ is a complexification of $K_0$. 

If $X$ is a subset of $\B$, $\overline{X}$ in this section will denote the closure of $X$ in $\B$ in the classical 
topology. If $X = K(v)$ or $X = B(w)$, where $v \in V$, $B$ is any Borel subgroup of $G$, and
$w \in W$, the closure of $X$ in the classical topology coincides with that in the Zariski topology.

\subsection{The Matsuki duality}\lb{subsec-Matsuli}
By 
Matsuki duality \cite{Matsuki, Uzawa-flow} \cite[$\S$6]{RS93}, for each $v\in V$, there exists a unique $G_0$-orbit $G_0(v)$
in $\B$ such that $G_0(v) \cap K(v)$ is a single $K_0$-orbit.
The map $v \mapsto G_0(v)$ gives an identification of the set
$V=\B/K$ of $K$-orbits in $\B$ with the set $\B/G_0$ of $G_0$-orbits in $\B$. 
For $v \in V$, let $K_0(v) = G_0(v) \cap K(v).$

The following \prref{pr-vvp} is proved in \cite[Corollary 1.4]{Uzawa-flow} \cite[Theorem 6.4.5]{RS93}.

\bpr{pr-vvp}
For $v, v^\prime \in V$, one has 
\[
v \leq v^\prime \hs \mbox{iff} \hs   G_0(v^\prime) \subset \overline{G_0(v)}
\hs \mbox{iff} \hs G_0(v) \cap K(v^\prime) \neq \emptyset.
\]
\epr

Let $v_0 \in V_0$ and fix $B \in K(v_0)$. Recall  the map
$q_B:  G \rightarrow  \B: g \mapsto B^g$ for $g \in G$.
For $v \in V$, let
\[
BvG_0 = q_B^{-1}(G_0(v)), \hs  BvK_0 = q_B^{-1}(K_0(v)),
\]
and recall that $BvK = q_B^{-1}(K(v))$. Let $B_K = K \cap B$. 

\ble{le-2-1}
1) $K = K_0B_K$ and $G_0K = G_0B_K$.

2) For any $v \in V$, one has $BvG_0K = \bigsqcup_{v' \in V, v' \geq v} Bv^\prime K$.
\ele

\begin{proof} 1) Let $U = \{g \in G: \sigma(g) = g\}$ and $H = B \cap \sigma(B)$. Then $U$ is a compact real form
of $G$ and $H$ a maximal torus of $G$. Let $A = \{h \in H: \sigma(h) = h^{-1}\}$ and let $N$ be the uniradical of $B$. Then one has the
Iwasawa decomposition $G = UAN$ of $G$, and $U$, $A$, and $N$ are all $\theta$-invariant. 
Let $k \in K$ and write $k = uan$ with $u \in U, a \in A$, and $n \in N$.
It follows from $\theta(k) = k$ that $u \in K_0$ and $an \in B_K$. Thus $K = K_0 B_K$.
It follows from $K_0 \subset G_0$ that  $G_0K = G_0B_K$.

2) Let $v \in V$. Then $B v G_0K$ is a union of $(B, K)$-double cosets, and for $v^\prime \in V$, 
\[
Bv^\prime K \subset BvG_0K \hs \mbox{iff} \hs 
Bv^\prime K \cap BvG_0 \neq \emptyset,
\]
 which, by \prref{pr-vvp}, is equivalent to $v \leq v^\prime$.
\end{proof}

\subsection{Proof of \thref{th-2-new}}\lb{subsec-proof-th-2-new}
Let $v_0 \in V_0$ and let $B \in K(v_0)$. Let $v \in V$ and $w \in W$. We must prove that
\[
G_0(v) \cap B(w) \neq \emptyset \hs \mbox{iff} \hs v \leq m(w)\sdot v_0.
\]
Let $BwB = q_B^{-1}(B(w)) \subset G$. Then, by definition,
 $G_0(v) \cap B(w) \neq \emptyset$ if and only if
$(BvG_0) \cap (BwB) \neq \emptyset$. Since $B_K \subset B$, this last
statement is equivalent to 
$(Bv G_0 B_K) \cap (BwB) \neq \emptyset$, which, by 1) of \leref{le-2-1}, is in turn equivalent to
$(BvG_0K) \cap (BwB) \neq \emptyset$. By 2) of \leref{le-2-1}, $G_0(v) \cap B(w) \neq \emptyset$ if and only if
$v \leq v^\prime$ for some $v^\prime \in V_w$, where $V_w = \{v^\prime \in V: \; K(v^\prime) \cap B(w) \neq \emptyset\}$.

By \prref{pr-mwv0}, $m(w)\sdot v_0$ is the unique maximal element in $V_w$ with respect to the Bruhat order on $V$.
Thus $G_0(v) \cap B(w) \neq \emptyset$ if and only if $v \leq m(w)\sdot v_0$.

This finishes the proof of \thref{th-2-new}.

\bco{co-closures-1} Let $v_0 \in V_0$ and let $B \in K(v_0)$. Then for any $v \in V$ and $w \in W$, 
the following are equivalent:

1) $G_0(v) \cap B(w) \neq \emptyset$;

2) $\overline{G_0(v)} \cap B(w) \neq \emptyset$;

3) $G_0(v) \cap \overline{B(w)} \neq \emptyset$;

4) $\overline{G_0(v)} \cap \overline{B(w)} \neq \emptyset$;

5) $K(v) \cap \overline{B(w)} \neq \emptyset$;

6) $K_0(v) \cap \overline{B(w)} \neq \emptyset$.
\eco

\begin{proof} Clearly 1) implies 4). Assume 4). Then there exist
$v_1 \in V$ and $w_1 \in W$ with $v_1 \geq v$ and $w_1 \leq w$, such that 
$G_0(v_1) \cap B(w_1) \neq \emptyset$. By \thref{th-2-new}, $v_1 \leq m(w_1)\sdot v_0$, and by \leref{le-y1-y},
$v \leq v_1 \leq m(w_1)\sdot v_0 \leq m(w) \sdot v_0.$
Thus 1) holds by \thref{th-2-new}. Hence 1) is equivalent to 4) and consequently also equivalent to 2) and 3). 

By \thref{th-2-new} and \leref{le-Wvv-combinatorial}, both 1) and 5) are equivalent to $v \leq m(w)\sdot v_0$.
Since $K=K_0B_K$ by 1) of \leref{le-2-1}, 5) and 6) are equivalent.
\end{proof}

\ble{le-y-v-K0}
Let $v_0 \in V_0$ and let $B \in K(v_0)$. Then for any $v \in V$ and $w \in W$, 
the following are equivalent:

1) $K(v) \cap B(w)  \neq \emptyset$;

2) $K_0(v) \cap B(w)  \neq \emptyset$.
\ele

\begin{proof}
It follows from $K = K_0B_K$ that 1) and 2) are equivalent.
\end{proof}

\sectionnew{Review on $K$-orbits in $\B$, II}\lb{sec-review-II}
Keeping the notation as in $\S$\ref{sec-review-I}, we now review the canonical definitions of more structures on the
set $V$ of $K$-orbits in $\B$, notably the cross action, the Springer map, root types, and reduced decompositions.
All the results in this section, except \prref{pr-v0v0-2} in $\S$\ref{subsec-closed}, can be found in \cite{RS90, RS93}.

\subsection{The canonical maximal torus $H_{\rm can}$ and the $W$-action on $H_{\rm can}$}\lb{subsec-Hcan}
Recall from $\S$\ref{subsec-C} that $\C$ is the variety of all pairs $(B, H)$, 
where $B \in \B$ and $H \subset B$ is a maximal torus of $G$.
Let $\widetilde{\C} = \{(B, H, h) \in \C \times G: (B, H) \in \C, \, h \in H\}$, and let
$G$ act on $\widetilde{\C}$  by 
\[
(B, H, h)^g = (B^g, H^g, g^{-1}hg), \hs (B, H, h) \in \widetilde{\C}, \; 
g \in G.
\]
 Let $H_{\rm can} = \widetilde{\C}/G$ 
be the set of $G$-orbits in $\widetilde{\C}$, and let $\widetilde{\C} \to H_{\rm can}: (B, H, h) \to[B, H, h]$
be the canonical projection. Then for every $(B, H) \in \C$, one has the bijective map 
\[
T_{B, H}: \; H \longrightarrow H_{\rm can}: \; \; h \longmapsto [B, H, h], \hs h \in H.
\]
Since $T_{B^\prime, H^\prime}^{-1} \circ T_{B, H}=T_{B, H}^{B^\prime, H^\prime}: H \to H^\prime$ for any
$(B, H), (B^\prime, H^\prime) \in \C$, where $T_{B, H}^{B^\prime, H^\prime}$ is the isomorphism of tori 
in \eqref{eq-T-BBHH},
$H_{\rm can}$ has a well-defined structure of a torus, such that 
$T_{B, H}: H \to H_{\rm can}$ is an isomorphism of
tori for every $(B, H) \in \C$. We will call $H_{\rm can}$ the {\it canonical maximal torus} of $G$.

The canonical Weyl group $W$ acts on $H_{\rm can}$ through the action of $W_H$ on $H$  via the identifications
$\eta_{B, H}^{-1}: W \to W_H$ and $T_{B, H}^{-1}: H_{\rm can} \to H$  for any $(B, H) \in \C$, and 
the action is
independent of the choice of $(B, H) \in \C$. If
$(B, H), (B^\prime, H) \in \C$, then
\begin{equation}\lb{eq-waction}
p(B^\prime, B) ([B, H, h]) = [B^\prime, H, h], \hs h \in H.
\end{equation}
One can also take \eqref{eq-waction} as the definition of the canonical action of $W$ on $H_{\rm can}$.

\subsection{The canonical root system  of $G$}\lb{subsec-Ga}
Let $X_{\rm can}$ be the character group of $H_{\rm can}$.
For $(B, H) \in \C$, let $X_H$ be the character group of $H$, and let 
$R_{B, H}: X_H \to X_{\rm can}$ be the isomorphism induced by $T_{B, H}^{-1}: H_{\rm can} \to H$. 
The {\it canonical sets of roots, positive roots, and simple roots} of $G$, are the subsets of $X_{\rm can}$, 
respectively defined by
\[
\Delta =R_{B, H}(\Delta_H), \hs \Delta^+ = R_{B, H}(\Delta^+_{B, H}), \hs \Gamma = R_{B, H}(\Gamma_{B, H}),
\]
where $(B, H) \in \C$,  and $\Delta_H \supset \Delta^+_{B, H} \supset \Gamma_{B, H}$ 
are the subsets of the character group $X_H$ of $H$
consisting, respectively, of the roots of $(G, H)$, the positive roots, and the simple roots determined by $B$.
It is easy to check that 
the sets $\Ga \subset \Delta^+ \subset \Delta \subset X_{\rm can}$ are independent of the choice of $(B, H) \in \C$.

For $\alpha \in \Delta$, choose any $(B, H) \in \C$, let $\alpha_H = R_{B, H}^{-1}(\alpha) \in \Delta_H$, let $s_{\alpha_H} \in W_H$
be the reflection defined by $\alpha_H$, and let $s_\alpha = \eta_{B, H}(s_{\alpha_H}) \in W$. Then
$s_\alpha$ is independent of the choice of $(B, H) \in \C$. Moreover, $S = \{s_\alpha: \alpha \in \Gamma\}$.

\subsection{The automorphism $\th$ on $H_{\rm can}$, $\Ga$, and $W$}\lb{subsec-th-Hcan}
In this subsection, let $\th$ be any automorphism of $G$, not necessarily of order $2$. Then one has the well-defined map
\[
\th: \; H_{\rm can} \longrightarrow H_{\rm can}: \; \; [B, H, h] \longmapsto [\th(B), \th(H), \th(h)], \hs (B, H, h) \in \widetilde{\C}.
\]
If 
$(B, H) \in \C$, then, by definition,
$\th = T_{\th(B), \th(H)} \circ \th|_H \circ T_{B, H}^{-1}$,
where $\th|_H: H \to \th(H)$. Thus $\th$ is an automorphism of $H_{\rm can}$. 

It is clear from the definition that if $\th$ is inner, i.e., if there exists $g_1 \in G$ such that
$\th(g) = g_1 g g_1^{-1}$ for $g \in G$, then $\th$ induces the identity automorphism of $H_{\rm can}$.

We will use the same letter to denote the induced 
action of $\th $ on the set $\Delta$ of canonical roots. 
It follows from the definitions that 
$\th(\Delta^+) = \Delta^+$ and $\th(\Gamma) = \Gamma$.

The following map, again denoted by $\th$, is easily seen to be well-defined:
\[
\th: \;  W \longrightarrow W: \; \; \; p(B_1, B_2) \longmapsto  p(\th(B_1), \th(B_2)), \hs B_1, B_2 \in \B,
\]
Choose any $(B, H) \in \C$, and let $\th_{W_H}: W_{H} \to W_{\th(H)}$ be the group isomorphism
induced by $\th|_{N_G(H)}: N_G(H) \to N_{G}(\th(H))$. It follows from definitions that
\[
\th=\eta_{\th(B), \th(H)} \circ \theta_{W_H} \circ \eta_{B, H}^{-1}: \; \; W \longrightarrow W.
\] 
Thus  $\th$ is an automorphism of $W$. It also follows that $\th(s_\alpha) = s_{\th(\alpha)}$ for all $\alpha \in \Delta$.

\subsection{The identification $\gamma: \C_\th/K \to V$}\lb{subsec-gamma} Let $\T^\th$ be the set of all $\th$-stable
maximal tori in $G$. The following \leref{le-T-th} is proved in  \cite[Proposition 1.2.1]{RS93} and
\cite[Corollary 4.4]{Springer85}.
See also \cite[1.4(b)]{RS90}.

\ble{le-T-th}
Every $B \in \B$ contains some $H \in \T^\th$, and if $H_1, H_2 \in \T^\th$ are both contained in $B$, then
$H_1 = k^{-1}H_2 k$ for some $k \in B \cap K$.
\ele

Let $\C_\th=\{(B, H)\in \C: H \in \T^\th\}$ and let $K$ act on $\C_\th$ by \eqref{eq-G-C}. It follows from \leref{le-T-th} that
the well-defined map
\[
\gamma: \; \C_\th/K \longrightarrow V=\B/K: \; \mbox{{\it K}-orbit of} \; (B, H) \; \mbox{in} \; \C \longmapsto 
\mbox{{\it K}-orbit of} \; B \; \mbox{in} \; \B,
\]
is a bijection (\cite[Proposition 1.2.1]{RS93} and 
\cite[Corollary 4.4]{Springer85}).

\subsection{The cross action of $W$ on $V$}\lb{subsec-WV-cross} We follow \cite[$\S$2]{RS90} and \cite[$\S$1.7]{RS93}.
See also \cite{adams-cloux}.
For any $(B, H) \in \C$ and $w \in W$, there is a unique
$(B^\prime, H) \in \C$ such that $w = p(B^\prime, B)$. Define
\begin{equation}\lb{eq-W-on-C}
W \times \C \longrightarrow \C: \; \; w \cdot (B, H) = (B^\prime, H), \hs (B, H), (B^\prime, H) \in \C, \, w = p(B^\prime, B).
\end{equation}
If we fix $(B, H) \in \C$, then under the identifications $\eta_{B, H}: W_H \to W$ and 
$C_{B, H}:  H\backslash G \to \C: Hg \mapsto (B^g, H^g)$ for $g \in G$, the map in 
\eqref{eq-W-on-C} becomes 
\[
W_H \times H\backslash G \longrightarrow H \backslash G: \; \; \; (nH, Hg) \longmapsto Hng, \hs n \in N_G(H), \, g \in G,
\]
which is a left action of $W_H$ on $H\backslash G$. Thus \eqref{eq-W-on-C} is indeed  a left action of $W$ on $\C$.
Moreover, the $W$-action commutes with the right action of $G$ on $\C$ given in \eqref{eq-G-C}.

By \eqref{eq-W-on-C},  $\C_\th \subset \C$ is $W$-invariant. Thus one has a well-defined action of
$W$ on $\C_\th/K$.
Identifying $V$ with $\C_\th/K$ via $\gamma$, one 
gets a left action of $W$ on $V$, which
is called the {\it cross action} of $W$ on $V$ and will be denoted by 
\[
W \times V \longrightarrow V: \; \; \; (w, v) \longmapsto w \sdot v, \hs w \in W, \, v \in V.
\]

\subsection{The Springer map $\phi: V \to W$}\lb{subsec-Springer-map}
The {\it Springer map} 
$\phi: V \to W$, introduced by Springer in \cite{Springer85}, is defined (see \cite[Remark 1.8]{RS90} and \cite[Proposition 1.7.1]{RS93})
by
\begin{equation}\lb{eq-phi-def}
\phi(v) = p(B, \theta(B)), \hs v \in V, \; B \in K(v).
\end{equation}
The element $\phi(v) \in W$ for $v \in V$ is an important invariant 
of the $K$-orbit $K(v) \subset \B$. Recall that $V_0$ consists of all $v \in V$ such that $K(v)$ is closed in $\B$. 
Let $1$ be the identity element of $W$. The following 
\prref{pr-two} is from \cite[Proposition 1.4.2]{RS93} and \cite[Proposition 2.5]{RS90}.

\bpr{pr-two}
1) For $v \in V$, $v \in V_0$ if and only if $\phi(v) = 1$.

2)  If $v, v^\prime \in V$ are such that $\phi(v) = \phi(v^\prime)$, then 
$v$ and $v^\prime$ are in the same $W$-orbit.
\epr

\subsection{Using standard pairs}
Let $(B, H)\in \C$ be a standard pair, i.e., 
$\th(B) = B$ and $\th(H) = H$.  Then for any $g \in G$, $(B^g, H^g)\in \C_\th$ if and only if
$g \th(g)^{-1} \in N_G(H)$, and for such a $g \in G$, 
\[
p(B^g, \,\th(B^g)) = p(B^g, \, B^{\th(g)}) = p(B^{g \th(g)^{-1}},  B) = \eta_{B, H}(g \th(g)^{-1}H) \in W.
\]
Letting ${\mathcal V}_H = \{g \in G: \, g \th(g)^{-1} \in N_G(H)\}$, 
one thus has the identification
\begin{equation}\lb{eq-VH-2}
H \backslash {\mathcal V}_H/K \longrightarrow V: \; \; \; HgK \longmapsto \mbox{the}\; \mbox{{\it K}-orbit in} \; \B \;
\mbox{through} \; B^g.
\end{equation}
Under the identification of $V$ with $H \backslash {\mathcal V}_H/K$ in \eqref{eq-VH-2} and that 
of $W$ with $W_H$ by $\eta_{B,H}$, the Springer map $\phi: V \to W$ becomes  
$\phi(HgK) = g \th(g)^{-1}H$ for $g \in {\mathcal V}_H$, and the action of $W$ on $V$ becomes 
$(nH) \sdot (HgK) = HngK$ for $n \in N_G(H)$ and $g \in {\mathcal V}_H$. 

In \cite{RS90, RS93, Springer85}, most of the structures on $V$ are introduced
and their properties proved by using standard pairs. For example, the 
following \leref{le-phi-W} immediately follows from the the identification
\eqref{eq-VH-2} (see also \cite[Lemma 2.1]{RS90}).

\ble{le-phi-W}
One has $\phi(w \sdot v) = w \phi(v) \th(w)^{-1}$ for any $w \in W$ and $v \in V$.
\ele

\subsection{The closed $K$-orbits in $\B$}\lb{subsec-closed}
Let $W^\th = \{w \in W: \th(w) = w\}$. By \leref{le-phi-W}, 
\begin{equation}\lb{eq-W-v0v0}
\{w \in W: \; w \sdot v_0 = v_0^\prime\} = \{w \in W^\th: \; w \sdot v_0 = v_0^\prime\}, \hs \forall \, v_0, v_0^\prime \in V_0.
\end{equation}
By 
\prref{pr-two},
$W^\th$ acts transitively on $V_0$. 

For the rest of this subsection, we assume that $K$ is connected.
 We will relate the 
sets in \eqref{eq-W-v0v0} for $v_0, v_0^\prime \in V_0$ with the canonical Weyl group $W_K$ of $K$.
\prref{pr-v0v0-2} will be used in $\S$\ref{subsec-V0}, where we determine $Y_{v_0}(v_0^\prime)$
for every $v_0, v_0^\prime \in V_0$.

\ble{le-BB-v0v0}
Let $v_0, v_0^\prime \in V_0$. Then for any $B \in K(v_0)$ and $B^\prime \in K(v_0^\prime)$, there exists
$H \in \T^\th$ such that $H \subset B \cap B^\prime$.
\ele

\begin{proof} 
Since $K(v_0)$ and $K(v_0^\prime)$ are closed, it follows
from \cite[Theorem 1.4.3]{RS93} (see also \cite[Corollary 6.6]{Springer85})
that $B$ and $B^\prime$ are $\th$-stable.  Thus, $B\cap B^\prime$
is $\th$-stable and contains a maximal torus of $G$. Hence, by
\cite[Theorem 7.5]{St-end}, $B \cap B^\prime$ contains
a $\th$-stable maximal torus of $G$.
\end{proof}

Let $v_0, v_0^\prime \in V_0$, and consider the restriction
$p: K(v_0^\prime) \times K(v_0)\rightarrow W$. By definition,
\[
p(K(v_0^\prime) \times K(v_0))=
\{w \in W: w=p(B^\prime, B) \, \mbox{for some} \; B^\prime \in K(v_0^\prime), \, B \in K(v_0)\}.
\]

\ble{le-BB-v0v0-W}
For any $v_0, v_0^\prime \in V_0$, one has
\[
p(K(v_0^\prime) \times K(v_0))=\{w \in W: w \sdot v_0 = v_0^\prime\}.
\]
\ele

\begin{proof} By the definition of the $W$-action,  $\{w \in W: w \sdot v_0 = v_0^\prime\}
\subset p(K(v_0^\prime) \times K(v_0))$. Suppose that $w = p(B^\prime, B)$ for some
$B \in K(v_0)$ and $B^\prime\in K(v_0^\prime)$. By \leref{le-BB-v0v0} and by the definition of the $W$-action on $V$,
$w \sdot v_0 = v_0^\prime$.
\end{proof}

Let $\B_K$ be the variety of all Borel subgroups of $K$ and let $K$ act on $\B_K$ by $(B_K)^k:=k^{-1} B_K k$ for
$B_K \in \B_K$ and $k \in K$.

\ble{le-flag-K}
For any $v_0 \in V_0$, the map
\[
\I_{v_0}: \; \; \; K(v_0) \longrightarrow \B_K: \; \; \; B \longmapsto B \cap K, \hs B \in K(v_0)
\]
is a $K$-equivariant isomorphism.
\ele

\begin{proof} By \cite[5.1]{R82}
(see also \cite[Page 113]{RS93}), $B \cap K \in \B_K$ for every $B \in K(v_0)$. 
Thus 
$\I_{v_0}$ is well-defined. It is clear that $\I_{v_0}$ is $K$-equivariant. 
To show that $\I_{v_0}$ is an isomorphism, we show that $\I_{v_0}$ is bijective and that
its inverse is an isomorphism 
from $\B_K$ to $K(v_0)$. 
Fix $B \in K(v_0)$ and identify $(B \cap K)\backslash K \cong \B_K$ via 
\[
(B \cap K)\backslash K \longrightarrow \B_K: \; \; \; (B\cap K) k \longmapsto (B\cap K)^k, \hs k \in K.
\]
 Consider the action map
$\eta: K \to K(v_0): k \mapsto B^k$, $k \in K$. Then the morphism $\tilde{\eta}: (B \cap K)\backslash K \to K(v_0)$ induced by $\eta$
is the inverse of $\I_{v_0}$. By \cite[Proposition 6.7 and Corollary 6.1]{Borel-book}, $\tilde{\eta}$ is an isomorphism if
$\eta$ is separable. 

Let $\fb, \fk,$ and $\fg$
be the Lie algebras of $B, K$, and $G$ respectively, and let $d\theta: \fg \to \fg$ be the 
differential of $\theta$. Since ${\rm char}({\bf k}) \neq 2$,  $\th$
is semisimple on $G$ and $B$ \cite[Section 5.4]{Spbook}. By \cite[Theorem 5.4.4(ii)]{Spbook}, 
$\fk=\fg^{d\theta} = \{x \in \fg: d\theta(x) = x\}$, and the Lie algebra of $B \cap K$ coincides with
$\fb \cap \fk$. Applying \cite[Proposition 6.12]{Borel-book} to the quotient morphism 
$G \to B\backslash G: g \to B^g, g \in G$, one sees that  $\eta$ is separable.
\end{proof}

Identify $W_K=(\B_K \times \B_K)/K$ for 
the diagonal action of $K$ on $\B_K \times \B_K$.
For $x \in W_K$, let $\O_K(x) \subset \B_K \times \B_K$
be the corresponding $K$-orbit in $\B_K \times \B_K$. Let $\leq_K$ denote the Bruhat order on $W_K$.

Let $v_0, v_0^\prime \in V_0$, and let
\[
\J_{v_0, v_0^\prime} = (\I_{v_0^\prime} \times \I_{v_0})^{-1}: \; \; \B_K \times \B_K \longrightarrow K(v_0^\prime) 
\times K(v_0).
\]
For $x \in W_K$, let $\O_{K, v_0, v_0^\prime}(x)$ be the single $K$-orbit in $K(v_0^\prime) \times K(v_0)$ given by
\begin{equation}\lb{eq-OK00}
\O_{K, v_0, v_0^\prime}(x) = \J_{v_0, v_0^\prime} (\O_K(x))\subset K(v_0^\prime) \times K(v_0).
\end{equation}
One then  has the well-defined map
\begin{equation}\lb{eq-I-v0v0}
I_{v_0, v_0^\prime}: \; \; \; W_K \longrightarrow W: \; \; \; x \longmapsto 
p\left(\O_{K, v_0, v_0^\prime}(x)\right),
\hs x \in W_K. 
\end{equation}
Let $1$ be the identity element in $W_K$.   By definition,
\[
\O_{K, v_0, v_0^\prime}(1)=\{(B^\prime, B) \in K(v_0^\prime) \times K(v_0): \; B^\prime \cap K = B \cap K\}.
\]

\bde{de-yv0v0}
For $v_0, v_0^\prime \in V_0$, let $y_{v_0, v_0^\prime} \in W^\th$ be  given by
\begin{equation}\lb{eq-yv0v0}
y_{v_0, v_0^\prime} = I_{v_0, v_0^\prime}(1) =p(B^\prime, B), \hs \forall (B^\prime, B) \in \O_{K, v_0, v_0^\prime}(1).
\end{equation}
\ede

\bpr{pr-v0v0-2}
Let $v_0, v_0^\prime \in V_0$. 

1) The map $I_{v_0, v_0^\prime}$ is a bijection from $W_K$ onto $\{w \in W^\th: \, w \sdot v_0 = v_0^\prime\}$;

2) If $x_1, x_2 \in W_K$ are such that
$x_1 \leq_K x_2$, then $I_{v_0, v_0^\prime}(x_1) \leq I_{v_0, v_0^\prime}(x_2)$;

3) $y_{v_0, v_0^\prime} \in W^\th$ is the unique minimal element in the set $\{w \in W^\th: \, w \sdot v_0 = v_0^\prime\}$ with respect 
to the Bruhat order on $W$.
\epr

\begin{proof}
1) By \leref{le-BB-v0v0-W} and \eqref{eq-W-v0v0}, the image of $I_{v_0, v_0^\prime}$ is $\{w \in W^\th: \, w \sdot v_0 = v_0^\prime\}$.
To show that $I_{v_0, v_0^\prime}$ is injective, 
assume that $B, B_1 \in K(v_0)$ and $B^\prime, B_1^\prime \in K(v_0^\prime)$ are such that 
$p(B^\prime, B) = p(B_1^\prime, B_1) \in W$. We must show that $(B^\prime, B)$ and $(B_1^\prime, B_1)$ are in the
same $K$-orbit in $K(v_0^\prime) \times K(v_0)$ for the diagonal action of $K$.
Without loss of generality, we may assume that $B_1 = B$, and we need to show that $B^\prime = (B_1^\prime)^k$ for some $k \in B \cap K$.
Let $H, H_1 \in \T^\th$ be such that $H \subset B \cap B^\prime$ and $H_1 \subset B \cap B_1^\prime$. By \leref{le-T-th},
there exists $k \in B \cap K$ such that $H = (H_1)^k$. Thus $H \subset B \cap B^\prime \cap (B_1^\prime)^k$. By the assumption,
$p(B^\prime, B) = p(B_1^\prime, B) = p((B_1^\prime)^k, B^k) = p((B_1^\prime)^k, B)$.
Thus $B^\prime =(B_1^\prime)^k$.

2) Suppose that  $x_1, x_2 \in W_K$ are such that $x_1 \leq_K x_2$. For $i = 1, 2$,  let 
$w_i = I_{v_0, v_0^\prime}(x_i) \in W$, so that $\O_{K, v_0, v_0^\prime}(x_i) \subset \O(w_i)$, where recall from  
$\S$\ref{subsec-Bruhat-W} that $\O(w_i)$ is the $G$-orbit in $\B \times \B$ corresponding to $w_i \in W$. 
For a subset $X$ of $\B_K \times \B_K$ (resp. of $\B \times \B$), let $\overline{X}^{\B_K \times \B_K}$
(resp. $\overline{X}^{\B \times \B}$) denote the Zariski closure of $X$ in $\B_K \times \B_K$ (resp. in $\B \times \B$).
Since $\J_{v_0, v_0^\prime}: \B_K \times \B_K \to \B \times \B$ is a morphism, one has 
\begin{align*}
\O_{K, v_0, v_0^\prime}(x_1) &= \J_{v_0, v_0^\prime} (\O_K(x_1)) \subset \J_{v_0, v_0^\prime} \left(\overline{\O_K(x_1)}^{\B_K \times \B_K}\right)
\subset \overline{\J_{v_0, v_0^\prime} (\O_K(x_1))}^{\B \times \B} \\&\subset \overline{\O(w_2)}^{\B \times \B}.
\end{align*}
 Thus $\O(w_1) \cap \overline{\O(w_2)}^{\B \times \B} \neq \emptyset$, 
and hence $w_1 \leq w_2$.

3) follows directly from 1) and 2).
\end{proof}

\bre{re-WK-v0} For $v_0 \in V_0$, let $W_K(v_0) = \{w \in W: w\sdot v_0 = v_0\}$. Then for $v_0, v_0^\prime \in V_0$, the set $\{w \in W^\th: \, w \sdot v_0 = v_0^\prime\}$ coincides with the coset $y_{v_0, v_0^\prime} W_K(v_0)$ in $W^\th$. It is easy to see
that $I_{v_0, v_0}: W_K \to W_K(v_0)$ is a group isomorphism (see, for example, \cite[Proposition 2.8]{RS90}).  Hence, by 3) of \prref{pr-v0v0-2}, every
coset in $W^\th/W_K(v_0)$ has a unique minimal length representative.
  In case $K$ is disconnected,
this last assertion is no longer true in the case when
 $G=PGL(4)$ and $K$ has connected
component of the identity equal to the image of $GL(2) \times GL(2)$
in $G$.
\ere

\subsection{The involution $\th_v: H_{\rm can} \to H_{\rm can}$}\lb{subsec-thv}
Let $v \in V$. Choose any $(B, H) \in \C_\th$ such that $B \in K(v)$, and  
let $\th|_H: H \to H$ be the restriction of $\th$ to $H$. Define
\[
\th_v :=T_{B, H} \circ \th|_H \circ T_{B, H}^{-1}: \;H_{\rm can} \longrightarrow H_{\rm can}: \; \; [B, H, h] \longmapsto [B, H, \th(h)], 
\hs h \in H.
\]
For another $(B^\prime, H^\prime) \in \C_\th$ such that $B^\prime \in K(v)$, there exists $k \in K$ such that
$B^\prime = B^k$ and $H^\prime = H^k$, and for any $h \in H$, $[B, H, h] = [B^k, H^k, k^{-1}hk] \in H_{\rm can}$, and 
\[
[B^k, H^k, \th(k^{-1}hk)] = [B^k, H^k, k^{-1}\th(h) k] = [B, H, \th(h)].
\]
Thus $\th_v: H_{\rm can} \to H_{\rm can}$ is  independent of 
the choice of $(B, H) \in \C_\th$. By \eqref{eq-waction},
\begin{equation}\lb{eq-thetav}
\th_v = \phi(v) \th: \;\; \; H_{\rm can} \longrightarrow H_{\rm can}.
\end{equation}
The induced involution on the
set $\Delta$ of canonical roots of $G$ (see $\S$\ref{subsec-Ga}) will also be denoted by $\th_v$.

\subsection{The subsets $p(s_\alpha, v)$ and root types}\lb{subsec-local}
Let $\alpha \in \Ga$. A parabolic subgroup of $G$ is said to be of {\it type} $\alpha$ if it is
of the form $B \cup Bs_\alpha B$ for some $B \in \B$. Let ${\mathcal P}_\alpha$ be the variety of all
parabolic subgroups of $G$ of type $\alpha$.
The action of $G$ on ${\mathcal P}_\alpha$ by conjugation is transitive, and one has the $G$-equivariant surjective
morphism
\[
\pi_\alpha: \; \; \B \longrightarrow {\mathcal P}_\alpha: \;  \; \pi_\alpha(B) = B \cup (B s_\alpha B), \hs B \in \B.
\]
For $v \in V$, let 
\[
p(s_\alpha, v) = \{v^\prime \in V: \; K(v^\prime) \subset \pi_\alpha^{-1}(\pi_\alpha(K(v))\}.
\]
It is well-known (see, for example, \cite[$\S$2.4]{RS93} and \cite{RS-comp}) that 
for each  $\alpha \in \Gamma$ and $v \in V$,
the subset
$p(s_\alpha, v)$ of $V$ has either one or two or three elements, depending on the {\it type} of $\alpha$ relative to $v$,
 and that $m(s_\alpha)\sdot v$ is the unique maximal element
in $p(s_\alpha, v)$ with respect to the Bruhat order on $V$. 

The case analysis of $p(s_\alpha, v)$ and the definition of the type of $\alpha$ relative to $v$ given in 
\cite[$\S$2.4]{RS93} and \cite{RS-comp} make use of a standard pair $(B_0, H_0) \in \C$, but the results are independent of the 
choice of $(B_0, H_0)$. Based on the results from \cite[$\S$2.4]{RS93} and \cite{RS-comp}, we give the equivalent definitions of root types in 
\deref{de-types} and summarize the results on $p(s_\alpha, v)$ from \cite[$\S$2.4]{RS93} and \cite{RS-comp}
in the  following \prref{pr-summary-6-RS93}.

\bde{de-types}
Let $v \in V$.
An $\alpha \in \Ga$ is said to be {\it imaginary} (resp. {\it real}, 
{\it complex}) for $v$ if
$\th_v (\alpha) = \alpha$ (resp. $\th_v(\alpha) =- \alpha, \;\th_v(\alpha) \neq \pm \alpha$). 
A simple imaginary root $\alpha$ for $v \in V$ is said to be {\it compact} if $m(s_\alpha)\sdot v = v$ and {\it non-compact}
if $m(s_\alpha)\sdot v \neq v$. A simple non-compact imaginary root $\alpha$ is said to be {\it cancellative}
if $s_\alpha \sdot v = v$ and {\it non-cancellative} if $s_\alpha \sdot v \neq v$ 
(see \cite[$\S$2.4]{Springer-invariants}). A simple real root 
$\alpha$ for $v \in V$ is said to be {\it cancellative} if $p(s_\alpha, v)$ has two elements and {\it non-cancellative}
if $p(s_\alpha, v)$ has three elements.  We will use the following notation.
\begin{align*}
I^c_v &= \{\alpha \in \Ga: \; \alpha \; \mbox{is compact imaginary for}\, v\},\\
I_v^{n,=}&= \{\alpha \in \Ga: \; \alpha \; \mbox{is non-compact imaginary and cancellative for}\, v\},\\
I_v^{n,\neq}&= \{\alpha \in \Ga: \; \alpha \; \mbox{is non-compact imaginary and non-cancellative for}\, v\},\\
R_v^{=}&= \{\alpha \in \Ga: \; \alpha \; \mbox{is real and cancellative for}\, v\},\\
R_v^{\neq}&= \{\alpha \in \Ga: \; \alpha \; \mbox{is real and non-cancellative for}\, v\},\\
C_v^+ &=\{\alpha \in \Ga: \; \alpha \; \mbox{is complex for}\, v \; \mbox{and} \;
\th_v(\alpha) \in \Delta^+\},\\
C_v^- &=\{\alpha \in \Ga: \; \alpha \; \mbox{is complex for}\, v \; \mbox{and} \;
\th_v(\alpha) \in -\Delta^+\}.
\end{align*}
We also set $I_v^n = I_v^{n, =} \cup I_v^{n, \neq}$, $I_v = I_v^c \cup I_v^n$, and $R_v = R_v^{=} \cup R_v^{\neq}$.
\ede

\bpr{pr-summary-6-RS93}
Let $v \in V$ and $\alpha \in \Gamma$.

{\it Case 1),} $\alpha \in I_v^c$. Then $p(s_\alpha, v) = \{v\}$. 

{\it Case 2),} $\alpha \in I_v^{n, =}$. Then $s_\alpha \sdot v
=v \neq m(s_\alpha)\sdot v$, and $p(s_\alpha, v) = \{v, \; m(s_\alpha) \sdot v\}$.
Moreover, $\alpha \in R_{m(s_\alpha)\cdot v}^{=}$.

{\it Case 3),} $\alpha \in I_v^{n, \neq}$. Then
$v$, $s_\alpha \sdot v$, and $m(s_\alpha)\sdot v$ are pair-wise distinct, and $p(s_\alpha, v) = \{v, \;
s_\alpha \sdot v, \;  m(s_\alpha) \sdot v\}$. Moreover, 
$\alpha \in I_{s_\alpha \cdot v}^{n, \neq}$ and $\alpha \in R_{m(s_\alpha)\cdot v}^{\neq}$.

{\it Case 4),} $\alpha \in R_v^{=}.$ Then there exists $v^\prime \in p(s_\alpha, v)$, $v^\prime \neq v$,
such that $s_\alpha \sdot v^\prime = v^\prime, \, v=m(s_\alpha) \sdot v^\prime =
s_\alpha \sdot v$, and
$p(s_\alpha, v) = \{v^\prime, \; v\}$. Moreover, $\alpha \in I_{v^\prime}^{n,=}$.

{\it Case 5),} $\alpha \in R_v^{\neq}$. Then there exists $v^\prime \in p(s_\alpha, v)$
such that $v^\prime, \, s_\alpha \sdot v^\prime$ and $v=m(s_\alpha) \sdot v^\prime=s_\alpha \sdot v$ 
are pair-wise distinct, and
$p(s_\alpha, v) = \{v^\prime, \, s_\alpha \sdot v^\prime, \, v\}$. Moreover, 
$\alpha \in I_{v^\prime}^{n, \neq} \cap  I_{s_\alpha\cdot v^\prime}^{n,\neq}$.

{\it Case 6),} $\alpha \in C_v^+.$  Then
$m(s_\alpha)\sdot v = s_\alpha \sdot v \neq v$, and $p(s_\alpha, v) = \{v, \, m(s_\alpha) \sdot v\}$. 
Moreover, $\alpha \in C^{-}_{m(s_\alpha)\cdot v}$.

{\it Case 7),} $\alpha \in C_v^-.$ Then
$m(s_\alpha) \sdot v = v \neq s_\alpha \sdot v$, and $p(s_\alpha, v) = \{v, \, s_\alpha \sdot v\}$.
Moreover, $\alpha \in C^{+}_{s_\alpha\cdot v}$.
\epr

\ble{le-122-RS93} \cite[Page 122]{RS93} Let $v, v^\prime \in V$ and $\alpha \in \Gamma$. Then
$v^\prime \in p(s_\alpha, v)$ if and only if $m(s_\alpha)\sdot v^\prime = m(s_\alpha)\sdot v$. 
Moreover, $p(s_\alpha, v^\prime) = p(s_\alpha, v)$ for all $
v^\prime\in p(s_\alpha, v).$
\ele

\ble{le-mv} \cite[$\S$3.2]{RS93} For $\alpha \in \Gamma$ and $v \in V$, $m(s_\alpha)\sdot v \neq v$ if and only if
$\alpha \in C_v^+ \cup I_v^n$. If $\alpha \in C_v^+$, then $\phi(m(s_\alpha) \sdot v) =s_\alpha \phi(v) s_{\theta(\alpha)}$,
and $l(\phi(m(s_\alpha) \sdot v)) =l(\phi(v))+2$.
If $\alpha \in I_v^n$, then $\phi(m(s_\alpha) \sdot v) = s_\alpha \phi(v) = 
\phi(v) s_{\theta(\alpha)}> \phi(v)$ and $\phi(s_\alpha \sdot v) = \phi(v)$.
\ele

\subsection{Reduced decompositions and subexpressions}\lb{subsec-reduced} 
We refer to \cite[$\S$5-7]{RS90} and \cite[$\S$4]{RS93} for more detail on this subsection.

\bde{de-deduced}\cite[Definition 3.2.3]{RS93} Let $v \in V$. A {\it reduced decomposition} of $v$ is a pair
$\bvs$, where ${\bf v} = (v_0, v_1, \ldots, v_k)$ is a sequence in $V$ and ${\bf s} = (s_{\alpha_1}, \ldots, s_{\alpha_k})$
is a sequence of simple reflections in $W$, such that $v_0 \in V_0$, $v_k = v$, and for each $j \in [1, k]$,
\[
\alpha_j \in C_{v_{j-1}}^+ \cup I_{v_{j-1}}^n \;\;\; \mbox{and} \;\;\; 
v_{j} = m(s_{\alpha_j}) \sdot v_{j-1}.
\]
The integer $k$ is called the length of the reduced decomposition
$\bvs$.
\ede

Every $v \in V$ has a reduced decomposition  and all reduced decompositions of $v$ have the same length, which
will be denoted by $l(v)$ and called the {\it length} of $v$ (see \cite[3.2 and Page 113]{RS93}).
\bde{de-sub}\cite[Definition 4.3]{RS93} Let $v \in V$ and let 
$({\bf v} = (v_0, v_1, \ldots, v_k), \, 
{\bf s} = (s_{\alpha_1}, \ldots, s_{\alpha_k}))$
be a reduced decomposition for $v$. A {\it subexpression} of $\bvs$ is a sequence ${\bf u} = (u_0, u_1, \ldots, u_k)$ in $V$ 
such that $u_0 = v_0$ and  one of the following 
holds for each $j \in [1, k]$: 

Case 1), $u_j = u_{j-1}$;

Case 2), $\alpha_j \in C_{u_{j-1}}^+ \cup I_{u_{j-1}}^n$ and $u_j = m(s_{\alpha_j})\sdot u_{j-1}$;

Case 3), $\alpha_j \in I_{u_{j-1}}^{n, \neq}$  and
$u_j = s_{\alpha_j}\sdot u_{j-1}$.

\noindent
In this case, $u_k$ is called the {\it final term} of the subexpression ${\bf u}= (u_0, u_1, \ldots, u_k)$.
\ede


\bpr{pr-sub-RS} \cite[Proposition 4.4]{RS93} Let $v, v^\prime \in V$ and let $\bvs$ be a reduced decomposition for 
$v$. Then $v^\prime \leq v$ if and only if there exists a subexpression  of $\bvs$ with final term $v^\prime$.
\epr

\sectionnew{Analysis on $\Yvv$ and proof of \thref{th-1-new}}\lb{sec-analysis-Yvv}

\subsection{The set $Y_{v_0}$}\lb{subsec-Yv0}
For $v_0 \in V_0$, let $Y_{v_0} = \bigcup_{v \in V} \Yvv$. 
Recall that $l$ denotes both the length function on $W$ and the length function on $V$.

\ble{le-Yvv-length}
Let $v_0 \in V_0$ and $y \in W$. Then $y \in Y_{v_0}$ if and only if $l(m(y)\sdot v_0) = l(y)$.
\ele

\begin{proof}
It is clear that $l(m(y)\sdot v_0) \leq l(y)$ for any $y\in W$.
Assume first that $y \in \Yvv$ for some $v \in V$.  If $l(m(y)\sdot v_0) < l(y)$, then
there exists $y_1 < y$ such that $m(y_1) \sdot v_0 = m(y) \sdot v_0$, so $y_1 \in \Wvv$, which is a contradiction.  
Thus $l(m(y)\sdot v_0) = l(y)$. Conversely, assume that $l(m(y)\sdot v_0) = l(y)$. Let $v = m(y) \sdot v_0$.
Then $y \in  \Yvv$.
\end{proof}

\bln{nota-bbvs}
Let $y \in W$. If $y = s_{\alpha_k} \cdots  s_{\alpha_1}$ is a reduced word of
$y$, the sequence ${\bf s} = (s_{\alpha_1}, \ldots, s_{\alpha_k})$ is also called a reduced word of $y^{-1}$.
Let 
\[
R(y^{-1})=\{{\bf s}: \;{\bf s} \; \mbox{is a reduced word of} \; y^{-1}\}.
\] 
Let $v_0 \in V_0$ and $y \in Y_{v_0}$. For   ${\bf s} = (s_{\alpha_1}, \ldots, s_{\alpha_k}) \in R(y^{-1})$, let
\begin{equation}\lb{eq-bvs-y}
{\bf v}_{v_0}({\bf s}) = (v_0, \; m(y_1)\sdot v_0,  \; \cdots, \; m(y_k) \sdot v_0),
\end{equation}
where $y_j =s_{\alpha_j} \cdots  s_{\alpha_1}$ for $j \in [1, k]$.
Then $\bbvs$ is a reduced decomposition for $m(y)\sdot v_0$, which will be called {\it 
the reduced decomposition of $m(y) \sdot v_0$ associated to ${\bf s}$}.
\eln

\begin{proof}
If $\bbvs$ is not a reduced decomposition for $m(y)\sdot v_0$, then $l(m(y) \sdot v_0) < l(y)$, and so $y \not\in Y_{v_0}$
by \leref{le-Yvv-length}.
\end{proof}

\subsection{Local analysis of $\Yvv$, Part I}\lb{subsec-key-lemma}
Recall the monoidal operation $\ast$ on $W$ in $\S$\ref{subsec-MWS}.

\ble{le-Min-0}
Let $v_0 \in V_0$, $v \in V$, and $\alpha \in \Gamma$.

1) If $w \in W_{v_0}( m(s_\alpha) \sdot v)$, then $w \in W_{v_0}(v^\prime)$ for all $v^\prime \in p(s_\alpha, v)$;

2) If $w \in \Wvv$, then $s_\alpha \ast w \in W_{v_0}(v^\prime)$ for all $v^\prime \in p(s_\alpha, v)$. 
\ele

\begin{proof} 1) follows from the fact that $v^\prime \leq m(s_\alpha) \sdot v$ for all $v^\prime \in p(s_\alpha, v)$.

2). Assume that $w \in \Wvv$. Then by \leref{le-y1-y},
\[
m(s_\alpha) \sdot v \leq m(s_\alpha) \sdot m(w) \sdot v_0 = m(s_\alpha \ast w) \sdot v_0.
\]
Thus $s_\alpha \ast w \in W_{v_0}(m(s_\alpha) \sdot v)$. 2) now follows from 1).
\end{proof}

The  following \leref{le-Min-key} is the key in proving \thref{th-1-new}.

\ble{le-Min-key} Let $v_0 \in V_0$, $v \in V$, and $y \in \Yvv$. Let $\alpha \in \Gamma$ and $y^\prime \in W$ be such that  
$y = s_\alpha y^\prime > y^\prime$. Then there exists $u \in p(s_\alpha, \, v)\backslash \{v,  \,m(s_\alpha)\sdot v\}$ such that
$y^\prime \in Y_{v_0}(u)$. Moreover,  there are three possibilities:

1) $\alpha \in C_u^+$ and  $v = m(s_\alpha) \sdot u$. In this case, $\alpha \in C_v^-$;

2) $\alpha \in I_v^n$ and $v = m(s_\alpha)\sdot u$. In this case, $\alpha \in R_v$;

3) $\alpha \in I_{u}^{n, \neq}$ and $v = s_\alpha \sdot u$. In this case, $\alpha \in I_v^{n, \neq}$. 
\ele

\begin{proof}
Let $y^\prime = s_{\alpha_{k-1}} \cdots s_{\alpha_1}$ be a reduced word of $y^\prime$. Then 
\[
{\bf s} = (s_{\alpha_1},  
\ldots, s_{\alpha_{k-1}}, s_{\alpha_k})\in R(y^{-1}),
\]
 where $\alpha_k = \alpha$.  Consider the reduced decomposition
$\bbvs$ of $m(y) \sdot v_0$.  By \prref{pr-sub-RS}, there is a subexpression
 ${\bf u}
=(u_0=v_0, \,u_1,,\ldots,,u_{k-1}, \, u_k)$  of $\bbvs$ with final term $u_k = v$, and we let $u = u_{k-1}$.
Since $y^\prime < y$ and $u \leq m(y^\prime)\sdot v_0$ by \prref{pr-sub-RS}, $u \neq v$ and $u \neq m(s_\alpha)\sdot u=m(s_\alpha)\sdot v$.
By the definition of a subexpression, one has either 1), 2), or 3). By \prref{pr-summary-6-RS93}, $\alpha \in C_v^-$ in 1),
$\alpha \in R_v$ in 2), and
$\alpha \in I_v^{n, \neq}$ in 3). 

It remains to show that $y^\prime \in Y_{v_0}(u)$. Since $u \leq m(y^\prime)\sdot v_0$, $y^\prime \in W_{v_0}(u)$.
Suppose that $y^\prime \notin
Y_{v_0}(u)$. Let $y^{\prime\prime} \in W_{v_0}(u)$ be such that $y^{\prime\prime} < y^\prime$. By \leref{le-Min-0}, 
$s_\alpha \ast y^{\prime\prime} \in \Wvv$. Since $s_\alpha \ast y^{\prime\prime} \leq s_\alpha \ast y^\prime = y$ and 
$l(s_\alpha \ast y^{\prime \prime}) \leq 1 + l(y^{\prime\prime}) < 1 + l(y^\prime) = l(y)$, we have 
$s_\alpha \ast y^{\prime\prime} < y$, which is a contradiction.
Thus $y^\prime \in Y_{v_0}(u)$.
\end{proof}

\subsection{Proof of \thref{th-1-new}}\lb{subsec-proof-th-2}
Assume that $v, v^\prime \in V$ are such that $\phi(v) = \phi(v^\prime)$ and  that there exists $y \in \Yvv
\cap Y_{v_0}(v^\prime)$. We use induction on $l(y)$ to show that $v = v^\prime$.

If $l(y) = 0$, then $y = 1$. By \leref{le-Wv-W}, $v=v^\prime =v_0$. 
Assume now that $l(y) \geq 1$, and choose $\alpha \in \Gamma$ such that
$y':= s_\alpha y < y$.
By \leref{le-Min-key}, there exist
$u \in p(s_\alpha, v)$ and $u^\prime \in p(s_\alpha, v^\prime)$ 
such that $y^\prime \in Y_{v_0}(u) \cap Y_{v_0}( u^\prime)$, and 
\[
\alpha \in \left(C_v^- \cup R_v \cup I_v^{n, \neq}\right) \cap 
\left(C_{v^\prime}^- \cup R_{v^\prime} \cup  I_{v^\prime}^{n, \neq}\right).
\]
Since $\phi(v) = \phi(v^\prime)$, one has by \eqref{eq-thetav} that
$\th_v(\alpha)  = \th_{v^\prime}(\alpha)$. Thus 
$\alpha \in C_v^-$ (resp. $R_v, \, I_v^{n, \neq}$) if and only if $\alpha \in C_{v^\prime}^-$ (resp. $R_{v^\prime}, \; I_{v^\prime}^{n, \neq}$).
We now look at the cases separately.

Case 1): $\alpha \in C_v^-$. In this case, $\alpha \in C_{v^\prime}^-$, $v = m(s_\alpha)\sdot u=s_\alpha \sdot u$, and 
$v^\prime =m(s_\alpha)\sdot u^\prime=
s_\alpha \sdot u^\prime$. By \leref{le-mv},
$\phi(u) = s_\alpha \phi(v) \theta(s_\alpha)
=s_\alpha \phi(v^\prime) \theta(s_\alpha)
=\phi(u^\prime).$
By the induction assumption, $u = u^\prime$. Thus $v = s_\alpha \sdot u = s_\alpha \sdot u^\prime = v^\prime$.

Case 2):  $\alpha \in R_v$. In this case, $\alpha \in R_{v^\prime}$, $v=m(s_\alpha)\sdot u$, and 
$v^\prime=m(s_\alpha) \sdot u^\prime$. By
\leref{le-mv}, $\phi(u) = s_\alpha \phi(v) = s_\alpha \phi(v^\prime) = \phi(u^\prime)$.
By the induction assumption, $u = u^\prime$. Thus $v = m(s_\alpha) \sdot u = m(s_\alpha) \sdot u^\prime = v^\prime$.

Case 3): $\alpha \in I_v^{n,\neq}$. In this case, $\alpha \in I_{v^\prime}^{n, \neq}$,
$u = s_\alpha \sdot v$, and $u^\prime = s_\alpha \sdot v^\prime$. By \leref{le-mv},
$\phi(u) = \phi(v) =\phi(v^\prime) = \phi(u^\prime)$. By the induction assumption, $u = u^\prime$.
Thus $v = s_\alpha \sdot u = s_\alpha \sdot u^\prime = v^\prime$.

This finishes the proof of \thref{th-1-new}.

\subsection{Local analysis on $\Yvv$, Part II}
The following \leref{le-Min-key-1} strengthens \leref{le-Min-key} and completes the local analysis on 
the sets $\Yvv$.

\ble{le-Min-key-1}
The element $u$ in \leref{le-Min-key} is unique.  If 2) of  \leref{le-Min-key} occurs and if $\alpha \in I_{u}^{n,\neq}$,
then $y \in Y_{v_0}(s_\alpha \sdot u)$. If 3) of \leref{le-Min-key} occurs, then $y \in Y_{v_0}(m(s_\alpha)\sdot v)$.
\ele

\begin{proof}
The only case where $u$ might not be unique is when 2) in  \leref{le-Min-key} occurs and when $\alpha \in I_{u}^{n,\neq}$.
Assume this is the case.  Since $\phi(u) = \phi(s_\alpha \sdot u)$,  by \thref{th-1-new}, $y^\prime$ can not be in both
$Y_{v_0}(u)$ and $Y_{v_0}(s_\alpha \sdot u)$, so the choice of $u$ is unique. Moreover, by 1) of \leref{le-Min-0},
$y \in W_{v_0}(s_\alpha \sdot u)$. Let 
$y_1 \in Y_{v_0}(s_\alpha \sdot u)$ be such that $y_1 \leq y$. 
By 2) of \leref{le-Min-0}, $s_\alpha \ast y_1 \in \Wvv$. Since
 $s_\alpha \ast y_1 \leq s_\alpha \ast y = y$, one has $s_\alpha \ast y_1 = y$. If
$s_\alpha y_1 > y_1$, then $s_\alpha y_1 = s_\alpha \ast y_1 =y$, so $y_1 = s_\alpha y = y^\prime$, contradicting the fact that
$y^\prime$ can not be in both
$Y_{v_0}(u)$ and $Y_{v_0}(s_\alpha \sdot u)$. Thus $s_\alpha y_1 < y_1$, and hence $y = s_\alpha \ast y_1=y_1 \in Y_{v_0}(s_\alpha \sdot u)$.

Assume now that 3) of \leref{le-Min-key} occurs. Then $y =s_\alpha y^\prime \in W_{v_0}(m(s_\alpha) \sdot v)$ by 2) of \leref{le-Min-0}. 
Let $y_2 \in Y_{v_0}(m(s_\alpha) \sdot v)$ be such that
$y_2 \leq y$.  By 1) of \leref{le-Min-0}, $y_2 \in \Wvv$. Thus $y = y_2 \in Y_{v_0}(m(s_\alpha) \sdot v)$. 
\end{proof}

\subsection{Elements of $\Yvv$ and subexpressions}
We now prove the following key property of $y \in \Yvv$ in terms of subexpressions of reduced decompositions of $m(y)\sdot v_0$.

\bpr{pr-unique}
Let $v_0 \in V_0, v \in V$ and $y \in \Yvv$. Then for any reduced word ${\bf s}=(s_{\alpha_1}, \ldots, s_{\alpha_k})$ 
of $y^{-1}$, 
there is exactly one subexpression ${\bf u}=(v_0, u_1, \ldots, u_k )$ of the reduced decomposition $\bbvs$ of $m(y)\sdot v_0$ 
that has final term $v$.
Moreover, for any $1\leq j \leq k$, $u_j \neq u_{j-1}$ and $s_{\alpha_j} s_{\alpha_{j-1}} \cdots s_{\alpha_1} \in 
Y_{v_0}(u_j)$.
\epr

\begin{proof} By \prref{pr-sub-RS}, there is a subexpression
 ${\bf u}=(u_0=v_0, u_1, \ldots, u_{k-1}, u_k)$  
of the reduced decomposition $\bbvs$ of $m(y)\sdot v_0$ with final term $v$. Letting $\alpha = \alpha_k$ and $y^\prime = 
s_{\alpha_{k-1}} \cdots s_{\alpha_1}$, one knows  from 
\leref{le-Min-key} and  \leref{le-Min-key-1}
that $y^\prime \in Y_{v_0}(u_{k-1})$ and that $u_{k-1}$ is uniquely determined by the quadruple $(v_0, v, y, \alpha)$.
Proceeding inductively, one proves  \prref{pr-unique}.
\end{proof}

\bnota{nota-bfu}
For $v_0 \in V_0$, $v \in V$, $y \in \Yvv$, and ${\bf s} \in R(y^{-1})$, the unique 
subexpression of the reduced decomposition 
$\bbvs$ of $m(y) \sdot v_0$ (see \lnref{nota-bbvs})
with final term $v$ will be denoted by  ${\bf u}_{v_{\scriptscriptstyle 0}, v}({\bf s})$. 
\enota

\sectionnew{Admissible paths and the sets $\Yvv$ and $\Zvv$}\lb{sec-Yvv-Zvv}

\subsection{Admissible paths}\lb{subsec-admissible}
\prref{pr-unique} leads naturally
to the following notion of admissible paths.

\bde{de-1-path} Fix $v_0 \in V_0$.

i) For  $v \in V$,  an {\it admissible path from $v_0$ to $v$} is a pair $\bvs$, where 
${\bf v}=(v_0, v_1, \ldots, v_k=v)$ is a sequence in $V$, and ${\bf s}=(s_{\alpha_1},s_{\alpha_2}, \ldots, s_{\alpha_k})$
is a sequence of simple reflections in $W$, and $k \geq 1$,
such that 
one of the following holds for each $j \in [1, k]$:

1) $\alpha_j \in C_{v_{j-1}}^+ \cup I_{v_{j-1}}^n$ and $v_j = m(s_{\alpha_j})\sdot v_{j-1}$;

2) $\alpha_j \in I_{v_{j-1}}^{n, \neq}$  and
$v_j = s_{\alpha_j}\sdot v_{j-1}$.

\noindent
The pair $({\bf v} = \{v_0\}, \, {\bf s} = \emptyset)$ is called the trivial admissible path from $v_0$ to $v_0$ with the understanding that
$k = 0$ in this case. The set of all admissible paths from $v_0$ to $v$ is denoted by $\calP(v_0, v)$.

ii) For $({\bf v}, {\bf s})\in \calP(v_0,v)$ as in i), the number $k$ is called the
{\it length} of $\bvs$ and denoted by $l({\bf v}, {\bf s})$. 
 We also set $y_0({\bf v}, {\bf s})=1\in W$, $y_j({\bf v}, {\bf s}) = s_{\alpha_j} s_{\alpha_{j-1}} \cdots s_{\alpha_1}$ for $j \in [1, k]$,
and $y({\bf v}, {\bf s}) = y_k\bvs$.

iii) For $v \in V$,  $\bvs \in \calP(v_0, v)$  is said to be a {\it shortest
admissible path} from $v_0$ to $v$
 if $l\bvs\leq l({\bf v}^\prime, {\bf s}^\prime)$ 
for every   $({\bf v}^\prime, {\bf s}^\prime) \in \Pv$. 
The length of a shortest admissible path from $v_0$ to $v$ will be denoted by $l_{v_0}(v)$.
\ede

\bre{re-both}
We are using the symbol $\bvs$ to denote both a reduced decomposition of an element in $V$ and an admissible path in $V$. 
This should cause no confusion as we will always use the modifiers   ``reduced decomposition" or ``admissible path" in front of
$\bvs$. Moreover, it is clear from the definitions that a reduced decomposition of $v \in V$ starting from $v_0 \in V_0$
is an admissible path from $v_0$ to $v$.
\ere

\ble{le-ad-ast}
Let $v_0 \in V_0$ and $v \in V$. Let $\bvs \in \Pv$, where ${\bf s} = (s_{1}, \dots s_{k})$. 
Then $v \leq m(s_{k} \ast \cdots \ast s_{1}) \sdot v_0$.
\ele

\begin{proof} Let ${\bf v} = (v_0, v_1 \ldots, v_{k-1}, v_k = v)$. Then 
\[
v \leq m(s_{k}) \sdot v_{{k-1}} \leq \cdots \leq m(s_{k}) m(s_{{k-1}}) \cdots m(s_{1}) \sdot v_0 = 
m(s_{k} \ast \cdots \ast s_{1}) \sdot v_0.
\]
\end{proof}

\subsection{Minimal admissible paths and elements in $\Yvv$}\lb{subsec-Yvv-admissible}
By \deref{de-sub}, 
for $v_0 \in V_0, v \in V, y \in \Yvv$, and ${\bf s} \in R(y^{-1})$, the pair
$({\bf u}_{v_0, v}({\bf s}), {\bf s})$  (see \notaref{nota-bfu}) is an admissible path from $v_0$ to $v$.
In this section, we give a characterization of such admissible paths.

\bde{de-minimal-admissible}
Let $v_0 \in V_0$ and $v \in V$.
For a sequence ${\bf s}$ of simple reflections in $W$, let 
\[
\calP(v_0, v, {\bf s}) =\{\bvsp \in \calP(v_0, v): \; {\bf s}^\prime \; \mbox{is a subsequence of} \; {\bf s}\}.
\]
A path $\bvs \in \calP(v_0, v)$ is said to be {\it minimal} if $\bvs$ is the only member of $\calP(v_0, v, {\bf s})$.
The set of all minimal admissible paths from $v_0$ to $v$ will be denoted by $\calP_{\rm min}(v_0, v)$.
\ede

\ble{le-minimal-paths}
Let $v_0 \in V_0$ and $v \in V$.

1) If $y \in \Yvv$ and ${\bf s} \in R(y^{-1})$, then $({\bf u}_{v_0, v}({\bf s}), {\bf s}) \in \calP_{\rm min}(v_0, v)$;

2) If $\bvs \in \calP_{\rm min}(v_0, v)$ and  $y = y\bvs$, then $y \in \Yvv$, ${\bf s} \in R(y^{-1})$, and
$\bvs =({\bf u}_{v_0, v}({\bf s}), {\bf s})$.
\ele

\begin{proof}
1) Let $y \in \Yvv$ and ${\bf s} \in R(y^{-1})$. Write ${\bf s} = (s_1, \ldots, s_k)$,  and 
let ${\bf s}^\prime = (s_{i_1}, \ldots, s_{i_p})$ be a subsequence of ${\bf s}$.
Suppose that $\bvsp \in \Pv$. By \leref{le-ad-ast}, $v \leq m(s_{i_p} \ast \cdots \ast s_{i_1}) \sdot v_0$.
Since $s_{i_p} \ast \cdots \ast s_{i_1} \leq y$, one must have $s_{i_p} \ast \cdots \ast s_{i_1}=y$, which is possible only when
${\bf s}^\prime = {\bf s}$. It follows that
${\bf v}^\prime$ is a subexpression of the reduced decomposition $({\bf v}_{v_0}({\bf s}), {\bf s})$ of $m(y)\sdot v_0$. 
By \prref{pr-unique}, 
${\bf v}^\prime ={\bf u}_{v_0, v}({\bf s})$. This proves that $({\bf u}_{v_0, v}({\bf s}), {\bf s}) \in \calP_{\rm min}(v_0, v)$.

2) Let $\bvs \in \calP_{\rm min}(v_0, v)$ and let $y = y\bvs$. Let ${\bf s} = (s_1, \ldots, s_k)$. By \leref{le-ad-ast},
$v \leq m(s_k \ast \cdots \ast s_1) \sdot v_0$. Suppose that $s_k \ast \cdots \ast s_1 \neq s_k \cdots s_1$. Then
there exists a proper subsequence ${\bf s}^\prime = (s_{i_1}, \ldots, s_{i_p})$ of ${\bf s}$ such that 
${\bf s}^\prime$ is a reduced word of $y^\prime =s_k \ast \cdots \ast s_1=s_{i_p} \cdots s_{i_1}$. By passing to a subsequence of
${\bf s}^\prime$ if necessary, we can assume that $y^\prime \in \Yvv$. Then $({\bf u}_{v_0, v} ({\bf s}^\prime), {\bf s}^\prime)$
is an admissible path from $v_0$ to $v$ that is different from $\bvs$,  contradicting the assumption on $\bvs$. Thus
$y=s_k \ast \cdots \ast s_1 = s_k \cdots s_1$, so ${\bf s}$ is a reduced word of $y^{-1}$. The same arguments show that
$y \in \Yvv$. By \prref{pr-unique}, $\bvs =({\bf u}_{v_0, v}({\bf s}), {\bf s})$.
\end{proof}

For $v_0 \in V_0$ and $v \in V$, let
\[
{\mathcal Y}_{v_0}(v) = \{(y, \,{\bf s}): \, y \in \Yvv, \, {\bf s} \in R(y^{-1})\}.
\]
The following \prref{pr-Yvv} follows immediately from \leref{le-minimal-paths}.

\bpr{pr-Yvv} For any $v_0 \in V_0$ and $v \in V$, the map
\begin{equation}\lb{eq-Omega}
\Omega: \; \; \; {\mathcal P}_{\rm min}(v_0, v) \longrightarrow {\mathcal Y}_{v_0}(v): \; \; \; \bvs \longmapsto (y\bvs, \, {\bf s})
\end{equation}
is bijective, with inverse given by 
\begin{equation}\lb{eq-Omega-inverse}
\Omega^{-1}:\; \; \; {\mathcal Y}_{v_0}(v) \longrightarrow {\mathcal P}_{\rm min}(v_0, v): \; \; \; (y, {\bf s}) \longmapsto ({\bf u}_{v_0, v}({\bf s}), {\bf s}).
\end{equation}
\epr

\bde{de-sub-y}
For $v_0 \in V_0, v \in V$, and $(y, {\bf s}) \in {\mathcal Y}_{v_0}(v)$, 
we will call  
$({\bf u}_{v_0, v}({\bf s}), {\bf s})$ the
{\it admissible path from $v_0$ to $v$ 
associated to  ${\bf s}$}. 
\ede

\bco{co-Yvv}
For any $v_0 \in V$ and $v \in V$, one has
\[
\Yvv = \{y\bvs: \bvs \in {\mathcal P}_{\rm min}(v_0, v)\}.
\]
\eco

\subsection{Shortest admissible paths and elements in $\Zvv$} Let $v_0 \in V_0$ and $v \in V$. 
Let ${\mathcal P}_{\rm short}(v_0, v)$ denote the set of all shortest admissible paths from $v_0$ 
to $v$ (see  \deref{de-1-path}). It is clear that ${\mathcal P}_{\rm short}(v_0, v)
\subset {\mathcal P}_{\rm min}(v_0, v)$. Let
\[
{\mathcal Z}_{v_0}(v) = \{(z, \,{\bf s}): \, z \in \Zvv, \, {\bf s} \in R(z^{-1})\}.
\]
Then ${\mathcal Z}_{v_0}(v) \subset {\mathcal Y}_{v_0}(v)$.

\bpr{pr-Zvv}
For any $v_0 \in V_0$ and $v \in V$, the map $\Omega$ in \eqref{eq-Omega} restricts to
a bijection between ${\mathcal P}_{\rm short}(v_0, v)$ and ${\mathcal Z}_{v_0}(v)$.
\epr

\begin{proof} Since $l\bvs = l(y\bvs)$ for all $\bvs \in {\mathcal P}_{\rm min}(v_0, v)$,
\prref{pr-Zvv} follows directly from \prref{pr-Yvv}.
\end{proof}

\bco{co-Zvv}
For any $v_0 \in V_0$ and $v \in V$, one has
\[
\Zvv = \{y\bvs: \bvs \in {\mathcal P}_{\rm short}(v_0, v)\}.
\]
\eco

\bex{ex-monoidal-orbit}
Let $v_0 \in V_0$, and let
$M(W, S) \sdot v_0 = \{m(w)\sdot v_0: w \in W\}$. 
Suppose that  $v \in M(W, S) \sdot v_0$. Then 
\begin{equation}\lb{eq-Zvv-right}
\Zvv =\{z \in W: \; m(z) \sdot v_0 = v, \, l(z) = l(v)\}.
\end{equation}
Indeed, if $w \in \Wvv$, then 
$l(v) \leq l(m(w)\sdot v_0) \leq l(w).$
Since $v \in M(W, S)\sdot v_0$, there exists $z_0 \in W$ such that $v = m(z_0)\sdot v_0$ and $l(z_0) = l(v)$.
Thus $l(v)=\min\{l(w): w \in \Wvv\}$. Let $Z_{v_0}^\prime(v)$ be the set of the right hand side on
\eqref{eq-Zvv-right}. Then $Z_{v_0}^\prime(v) \subset \Zvv$. Conversely, if $z \in \Zvv$, then $l(z) = l(v)$, and it follows from 
$l(v)\leq l(m(z)\sdot v_0) \leq l(z)$ that $v = m(z)\sdot v_0$. Thus $\Zvv \subset Z_{v_0}^\prime(v)$. Hence
$\Zvv = Z_{v_0}^\prime(v)$. 
\eex

We will see in $\S$\ref{subsec-an-ex} that it can happen that $v \in M(W, S) \sdot v_0$ and $\Zvv$ is a proper subset of $\Yvv$.


\sectionnew{Examples}\lb{sec-examples}

\subsection{The case when $v \in V_0$}\lb{subsec-V0}
In this subsection, we assume that $K$ is connected. Fix $v_0, v_0^\prime \in V_0$. 
We will consider the set $Y_{v_0}(v_0^\prime)$.  

Recall from 
$\S$\ref{subsec-Wvvp} that $Y_{v_0}(v_0^\prime) = \min (W_{v_0}^\prime(v_0^\prime))$, where
\begin{equation}\lb{eq-w0w0p}
W_{v_0}^\prime(v_0^\prime)=\{w \in W: \; w=p(B^\prime, B) \; \mbox{for some} \; B^\prime \in K(v_0^\prime), \, B \in K(v_0)\},
\end{equation}
and 
$p: \B \times \B \to W=(\B \times \B)/G$ is the natural projection.
By \leref{le-BB-v0v0-W}, $W_{v_0}^\prime(v_0^\prime) = \{w \in W^\th: w \sdot v_0 = v_0^\prime\}.$
Recall from \deref{de-yv0v0} the element $y_{v_0, v_0^\prime} \in W^\th$.
By 3) of \prref{pr-v0v0-2}, $y_{v_0, v_0^\prime}$ is the unique minimal element in $W_{v_0}^\prime(v_0^\prime)$. Thus
$Y_{v_0}(v_0^\prime) = \{ y_{v_0, v_0^\prime} \}$
has one element.

\subsection{The case when there is a unique closed orbit}\lb{subsec-unique-closed}
Assume that $V_0=\{v_0\}$ has only one element. Then $M(W, S)\sdot v_0  = V$. In this case, for every $v \in V$,
one has
\[
\Yvv = \Zvv = \{z \in W: \; m(z) \sdot v_0 = v, \, l(z) = l(v)\}.
\]
Indeed, by \cite[Corollary 9.15]{RS90}, the Springer map $\phi: V \to W$ is injective. For any $v \in V$, since
$\phi(s_\alpha \sdot v) = \phi(v)$ for any $\alpha \in I_v^n$, one has
 $I_v^{n, \neq} = \emptyset$, and thus for any
admissible path $\bvs$ from $v_0$ to $v$, 2) in \deref{de-1-path} does not occur. Therefore if $y \in \Yvv$, then
$v = m(y)\sdot v_0$ and $l(y) = l(v)$, and so $y \in \Zvv$.

\bre{re-Springer-invariants} For any $v \in V$, and not assuming that $V_0$ has only one element, T. A. Springer
has studied in \cite{Springer-invariants} the set $\cup_{v_0 \in V_0}\{z \in W: m(z)\sdot v_0 = v, \, l(z) = l(v)\}$ as an invariant for $v$.
\ere

\subsection{The Hermitian symmetric case}\lb{subsec-hermitian}
 Assume that $G$ is simple and simply connected. By \cite[Definition 5.1.1]{RS93},
$(G, \th)$ is said to be of Hermitian symmetric type if the center of $K$ has positive dimension.

Assume that $(G, \th)$ is of Hermitian symmetric type.
By \cite[Theorem 5.12]{RS93}, there exists a standard pair $(B, H) \in \C$ and a parabolic subgroup $P$ of
$G$ containing $B$ such that $K$ is the unique Levi subgroup of $P$ containing $H$. Moreover, 
there exists $\alpha_0 \in \Gamma$ such that 
$P$ is of type $J = \Gamma\backslash \{\alpha_0\}$. 
Let $v_0 \in V$ such that $B \in K(v_0)$. Then every $\alpha \in J$ is compact imaginary for $v_0$.
We now study the set $\Yvv$ for every $v \in V$.

Let $W_J$ be the subgroup of $W$ generated by simple reflections corresponding to roots in $J$, and let
$W^J$ be the set of minimal length representatives for $W/W_J$ in $W$. Then $W^J$ parametrizes
the set of $(B, P)$-double cosets in $G$ via 
\[
W^J \ni d \longmapsto BdP:=B \eta_{B, H}^{-1}(d) P \subset G,
\]
where $\eta_{B, H}: W_H \to W$ is given in \eqref{eq-etaBH}.
Since every $(B, K)$-double coset in $G$ is contained in a 
unique $(B, P)$-double coset, we have the well-defined surjective map
\[
\nu: \; V \longrightarrow W^J: \; \; \; \nu(v) = d \in W^J \hs \mbox{{\rm if}}  \hs B v K \subset B d P.
\]
It is proved in \cite[Theorem 5.2.5]{RS93} that the map 
\[
\eta: \; V \longrightarrow W \times W^J: \; \; \eta(v) = (\phi(v), \; \nu(v))
\]
is injective.

\bpr{pr-hermitian-case}
Let $(G, \th)$ be of Hermitian symmetric type and let the notation be as above. Then for any $v \in V$, 
$\nu(v) \in W^J$ is the unique  element in $\Yvv$. 
\epr

\begin{proof} Note that $P = BK$. Let $v \in V$.
By the definition of
$\nu(v)$, one has
\[
BvK \subset B\nu(v) P \subset \overline{B \nu(v)P} = \overline{B \nu(v)BK}  =
 \overline{B\nu(v)B} \; K = \overline{B (m(\nu(v))\sdot v_0) K}.
\]
Thus $\nu(v) \in \Wvv$. 
Conversely, assume that $w \in \Wvv$. 
Write $w = d x$, where $d \in W^J$ and $x \in W_J$. Since
every $\alpha \in J$ is compact imaginary for $v_0$, $m(x) \sdot v_0 = v_0$. Thus 
$v \leq m(w)\sdot v_0 = m(d) \sdot v_0$, i.e., 
\[
BvK \subset \overline{B (m(d)\sdot v_0) K}=\overline{B d B} \; K = \overline{BdP}.
\]
Since $B v K \subset B \nu(v)P$, we have
$B \nu(v)P \cap \overline{BdP} \neq \emptyset$. Thus $\nu(v) \leq d \leq w$. This shows that
$\nu(v)$ is the unique minimal element in $\Wvv$.
\end{proof} 

Thus, our \thref{th-1-new} generalizes \cite[Theorem 5.2.5]{RS93}.

\subsection{An example where $\Zvv \neq \Yvv$}\lb{subsec-an-ex}

Let $G = SL(4, \Cset)$ and let $\theta \in {\rm Aut}(G)$ be given by $\th(g) = I_{2, 2} (g^t)^{-1} I_{2, 2}$, where
for $g \in G$, $g^t$ denotes the transpose of $g$, and $I_{2, 2} = {\rm diag}(I_2, -I_2)$ with $I_2$ being the
$2 \times 2$ identity matrix. Then $K = S(GL(2, \Cset) \times GL(2, \Cset))$. Using the ``kgb" command in the {\it Atlas of Lie groups} (www.liegroups.org) 
for the real form $SU(2, 2)$ of $G$, 
one knows that there are $6$ elements in $V_0$. Take $v_0 \in V_0$ to be the orbit labeled by $3$ in the Atlas
and let $v$ be the orbit labeled by $19$. With respect to $v_0$, 
 the simple roots $\alpha_1$ and $\alpha_3$ are noncompact, while
$\alpha_2$ is compact, where we use the Bourbaki labeling of roots.  Then $v = m(s_3s_2s_1) \sdot v_0 \in M(W, S) \sdot v_0$.
By \exref{ex-monoidal-orbit},  $l(z) = 3$ for every $z \in \Zvv$. 

On the other hand, let $y = s_2s_1s_2s_3 \in W$. Then $m(y) \sdot v_0 = v_{{\rm max}}$, where $v_{{\rm max}}$ is maximal element
in $V$. Since $v \leq v_{{\rm max}}$, $y \in \Wvv$. We claim that $y \in \Yvv$. Indeed, if $y \notin \Yvv$, then
there exists $y^\prime \in \Yvv$ with $y^\prime < y$. Since $l(y) \geq l(v) = 3$, we must have $l(y) = 3$. Now there
are exactly three subwords $y^\prime$ of the word $y= s_2s_1s_2s_3$ with length $3$, namely, $y^\prime = s_1s_2s_3$ or $s_2s_1s_3$, or 
$s_2s_1s_2$, and one  checks directly that neither of these three choices gives $v \leq m(y^\prime) \sdot v_0$. Thus
$y \in \Yvv$. Since $l(y) = 4$, $y \notin \Zvv$. We thus have an example
in which $\Zvv$ is a proper subset of $\Yvv$.


\end{document}